\documentclass[11pt, a4paper]{amsart}
\usepackage{amsfonts, amsmath, amssymb,latexsym,amsthm, mathrsfs, enumerate, tikz-cd}
\usepackage[utf8]{inputenc}
\usepackage{csquotes}
\usepackage[english]{babel}
\usepackage{epsfig}

\usepackage[margin=1.2in]{geometry}
\usepackage{graphicx}
\usepackage{listings}

\usepackage{hyperref}
\usepackage{xcolor}
\usepackage[capitalise]{cleveref}
\usepackage{mathtools}
\usepackage{thmtools}

\usepackage[
    backend=biber,
    style=alphabetic,
    url=true,
]{biblatex}
\newtheorem{thm}{Theorem}[section]
\newtheorem*{thm*}{Theorem}
\newtheorem*{prop*}{Proposition}
\newtheorem*{cor*}{Corol·lari}
\newtheorem*{conj*}{Conjecture}
\newtheorem*{ques*}{Question}
\newtheorem{prop}[thm]{Proposition}
\newtheorem{cor}[thm]{Corollary}
\newtheorem{lem}[thm]{Lemma}
\newtheorem{defn}[thm]{Definition}

\newtheorem{conj}{Conjecture}

\theoremstyle{definition}

\newtheorem{rmk}[thm]{Remark}

\newcommand{\Z}{\mathbb{Z}}
\newcommand{\R}{\mathbb{R}}
\newcommand{\QQ}{\mathbb{Q}}
\newcommand{\RP}{\mathbb{R}P}
\newcommand{\C}{\mathbb{C}}
\newcommand{\CP}{\mathbb{C}P}

\renewcommand{\ker}{\operatorname{Ker}}

\newcommand{\Diff}{\operatorname{Diff}}

\newcommand{\rank}{\operatorname{rank}}
\newcommand{\Aut}{\operatorname{Aut}}

\newcommand{\SO}{\operatorname{SO}}
\newcommand{\Spin}{\operatorname{Spin}}

\newcommand{\SU}{\operatorname{SU}}
\newcommand{\Sp}{\operatorname{Sp}}

\newcommand{\Zc}{\mathcal{Z}}

\newcommand{\A}{\operatorname{tor-sym}}
\newcommand{\Q}{\operatorname{q-sym}}
\newcommand{\Qr}{\operatorname{q-\rank}}

\newcommand{\acts}{\curvearrowright}
\newcommand{\core}{\operatorname{core}}

\usepackage{fancyhdr}

\author{Jordi Daura Serrano}
\address{Department of Mathematics and Statistics, University of Cyprus (UCY), Nicosia 1678 (Cyprus)}
\email{daura-serrano.jordi@ucy.ac.cy}
\author{Jordi Garriga Puig}
\address{Departament de Matemàtiques i Informàtica, Universitat de Barcelona (UB), Gran Via de les Corts Catalanes 585, 08007 Barcelona (Spain)}
\email{jordi.garrigapuig@ub.edu}
\author{Ioannis Gkeneralis}
\address{Department of Mathematics, Aristotle University of Thessaloniki (AUTh), Thessaloniki 54124 (Greece)}
\email{igkeneralis@math.auth.gr}

\subjclass[2020]{Primary 57S15; Secondary 57S25, 53C30.}

\keywords{quaternionic torus actions, quaternionic degree of symmetry, quaternionic rank, almost-effective actions, transformation groups, cohomogeneity one, maximal symmetry rank.}

\title{Quaternionic Torus Actions on Manifolds}
\begin{document}

\pagestyle{plain}
\lhead{\textbf{}}
\lfoot{}
\rhead{}
\cfoot{}
\rfoot{\thepage}

\begin{abstract}
We study smooth actions of quaternionic tori $\SU(2)^r$ on closed, connected smooth manifolds, establishing sharp bounds for $r$ in terms of the dimension of the manifold when the action is either effective, effective almost-free, or almost-effective. In particular, we prove that an almost-effective $\SU(2)^r$-action on an $n$-dimensional manifold satisfies $r\leq \lfloor2n/3\rfloor$, and we classify up to diffeomorphism all simply connected manifolds attaining the maximal possible rank. This gives a quaternionic counterpart of the maximal symmetry rank conjecture without any curvature assumption. We also establish the bounds $r\leq n/3$ for effective and effective almost-free quaternionic torus actions, and we compute the exact value for spheres, quoric manifolds, or manifolds of the form $(S^3)^a\times (S^1)^b$.
\end{abstract}
\maketitle

\pagenumbering{arabic} 
\setcounter{page}{1}

\section{Introduction}

The study of torus actions on manifolds is a classical subject in the theory of compact transformation groups, with applications across algebraic topology, differential geometry, and algebraic geometry. A common theme is determining how the topology and geometry of a manifold constrain the rank of tori that can act effectively on it. For example, a classical theorem states that if $T^r$ acts effectively on a closed, connected, smooth $n$-dimensional manifold, then $r\leq n$ and the bound is attained if and only if $M$ is diffeomorphic to $T^n$. In differential geometry, Grove's symmetry program studies isometric group actions of large compact Lie groups on manifolds of non-negative sectional curvature. In particular, the maximal symmetry rank conjecture predicts that if $T^r$ acts effectively and isometrically on a closed, simply connected, non-negatively curved $n$-manifold, then $r\leq\lfloor\tfrac{2n}{3}\rfloor$, with strong restrictions on the manifolds for which the equality is attained; see \cref{conj: maximal symmetry rank conjecture} for a precise formulation and \cite{grove2002geometry, searle2023symmetries} for surveys of the symmetry program. Toric topology provides another source of examples. Davis and Januszkiewicz \cite{davis1991convex} introduced quasitoric manifolds as topological analogues of smooth projective toric varieties. They are closed, smooth $2n$-manifolds with an effective $T^n$-action that is locally the standard $T^n$ action on $\C^n$ and such that the orbit space is a simple polytope. These and many other classes of manifolds connecting torus actions and combinatorics are presented in \cite{BuchstaberPanov2015}.

Actions of compact connected non-abelian Lie groups have also been extensively studied. For example, \cite{burghelea1975semisimple, saito1984semisimple, watabe1985semisimple, xu2003degree, xu2007degree} obtained bounds for the dimension of semisimple compact Lie groups acting effectively on some classes of manifolds, while Lawson and Yau \cite{lawson1974scalar} proved that a closed, connected, smooth manifold admitting an effective action of such a group admits a metric of positive scalar curvature. In toric topology, Hopkinson introduced quoric manifolds as quaternionic analogues of quasitoric manifolds in \cite{Hopkinsonthesis}. They have dimension $4n$ and carry locally regular actions of $\SU(2)^n$ whose orbit spaces are $n$-dimensional simple polytopes. Related constructions connecting actions of $\SU(2)^r$ with combinatorial objects can be found in \cite{scott1995quaternionic, gentiligorisarfatti, batakidisgkeneralistopology, batakidis2026tetraplectic}. We also note that the orbit space of actions of $\SU(2)^r$ on products of spheres is studied in \cite{ayzenberg2024onactions}.

These developments motivate our study of actions of $\SU(2)^r$ on manifolds, with the goal of examining how they resemble their abelian counterparts and which features of semisimple group actions they retain. We write 
$$Q:=\Sp(1) \cong \SU(2).$$
Since $Q$ is diffeomorphic to $S^{3}$, we call $Q^r$ a \textit{quaternionic torus} of rank $r$. 

Throughout the paper, all manifolds are closed, connected and smooth, and all actions are smooth. We denote a smooth action of a Lie group $G$ on a manifold $M$, $\Phi: G \longrightarrow \Diff(M)$, by $G\acts M$. The \textit{kernel} of the action is the subgroup $\ker(\Phi) \subseteq G$. Recall that an action is \textit{effective} if its kernel is trivial, \textit{almost-effective} if its kernel is finite, and \textit{almost-free} if all isotropy subgroups are finite. 

We introduce three invariants that determine the maximal $r$ such that $Q^r$ acts on a manifold $M$ under different assumptions on the actions.

\begin{defn}
\label{def: quaternionic invariants}
Let $M$ be a closed, connected smooth manifold. Then:
\begin{itemize}
    \item[1.] The \textbf{quaternionic degree of symmetry} of $M$ is 
    $$\Q(M) = \max \{ r : Q^{r} \acts M \textnormal{ effectively} \}.$$ 
    \item[2.] The \textbf{quaternionic rank} of $M$ is
    $$\Qr(M) = \max\{r:Q^r\acts M\textnormal{ effectively and almost-freely}\}.$$
    \item[3.] The \textbf{almost-effective quaternionic degree of symmetry} of $M$ is 
    $$\Q_{ae}(M)=\max\{r: Q^r \acts M \textnormal{ almost-effectively} \}.$$  
\end{itemize}
\end{defn}

The first two invariants are quaternionic counterparts of the toral degree of symmetry $\A(M)$ and the toral rank $\rank(M)$, defined using effective and almost-free actions of $T^r$, respectively (see \cite{hsiang2012cohomology, lee2010seifert, mundet2024discrete} and \cite{allday1993cohomological, felix2008algebraic} for background on the toral degree of symmetry and the toral rank, respectively). 
If $G$ is a compact, connected Lie group acting almost-effectively on $M$, then $G/\Gamma$ acts effectively on $M$, where $\Gamma$ denotes the kernel of the $G$-action on $M$. Thus, if $T^r$ acts almost-effectively on $M$ with kernel a finite group $\Gamma \subseteq T^{r}$, then $T^r/\Gamma\cong T^r$ acts effectively on $M$. Note that if $Q^r$ acts almost-effectively on a manifold $M$, then the restriction of the action to a maximal torus $T^r\leq Q^r$ is also almost-effectively, and therefore $\Q_{ae}(M)\leq \A(M)$. Hence, for a closed, connected smooth manifold $M$, we have the string of
inequalities $$\Qr(M)\leq \Q(M)\leq \Q_{ae}(M)\leq \A(M).$$

Our main result gives a sharp upper bound for the almost-effective quaternionic degree of symmetry and classifies, up to diffeomorphism, the simply connected manifolds for which the bound is attained.

\begin{thm}
\label{bound almost-effective quaternionic degree of symmetry}
Let $M$ be a closed, connected smooth $n$-dimensional manifold. Then $\Q_{ae}(M)\leq\tfrac{2n}{3}$. Moreover, if $M$ is simply connected and $\Q_{ae}(M)=\lfloor\tfrac{2n}{3}\rfloor$, then:
\begin{itemize}
    \item[1.] If $n=3k$, then $M$ is diffeomorphic to $(S^3)^k$.
    \item[2.] If $n=3k+1$, then $M$ is diffeomorphic to either
    \begin{itemize}
        \item[\textnormal{(a)}] $(S^3)^{k-1}\times S^2\times S^2$,  
        \item[\textnormal{(b)}] $(S^3)^{k-1}\times S^4$, or 
        \item[\textnormal{(c)}] $(S^3)^{k-2}\times S^7$.
    \end{itemize}
    \item[3.] If $n=3k+2$, then $M$ is diffeomorphic to $(S^3)^{k}\times S^2$.
\end{itemize}
\end{thm}

Every manifold in this classification admits an almost-effective quaternionic torus action that attains the bound and preserves the standard product metric, which has non-negative sectional curvature. This relates our theorem to the maximal symmetry rank conjecture, a central problem in the study of isometric torus actions on non-negatively curved manifolds.

\begin{conj}(Maximal symmetry rank conjecture, \cite{galaz2011low, escher2021torus, dong2024almost})
\label{conj: maximal symmetry rank conjecture}
Let $T^r$ act isometrically and effectively on a closed, simply connected, non-negatively curved $n$-dimensional Riemannian manifold $M$. Then $r\leq \tfrac{2n}{3}$. Moreover, if $r=\lfloor\tfrac{2n}{3}\rfloor$, then:
\begin{itemize}
    \item[1.] If $n=3k$, then $M$ is diffeomorphic to $(S^3)^k$.
    \item[2.] If $n=3k+1$, then $M$ is diffeomorphic to either
    \begin{itemize}
        \item[\textnormal{(a)}] $(S^3)^{k-1} \times S^4$,  
        \item[\textnormal{(b)}] $(S^3)^{k-2} \times S^7$, 
        \item[\textnormal{(c)}] $(S^3)^{k-3} \times S^5\times S^5$,
        \item[\textnormal{(d)}] $(S^3)^{k+1}/T^2$, or
        \item[\textnormal{(e)}] $N/T^1$, where $N$ is a manifold belonging to item $3$.
    \end{itemize}
    \item[3.] If $n=3k+2$, then $M$ is diffeomorphic to either
    \begin{itemize}
        \item[\textnormal{(a)}] $(S^3)^{k-1} \times S^5$, or
        \item[\textnormal{(b)}] $(S^3)^{k+1}/T^1$.
    \end{itemize}
\end{itemize} 
The torus actions defining the quotients are free and linear.
\end{conj}

Note that since the quotient of a free linear action of $T^1$ on $S^3$ is $S^3/T^1\cong S^2$, the cases $(S^3)^{k-1}\times S^2\times S^2$ and $(S^3)^k\times S^2$ correspond to the cases 2.(d) and 3.(b) respectively. \Cref{bound almost-effective quaternionic degree of symmetry} proves a quaternionic version of the maximal symmetry rank conjecture without any curvature assumption. Note that if $Q^r$ acts almost-effectively and isometrically on $M$, then we can construct an effective isometric action of $T^r$ on $M$. Thus, \cref{bound almost-effective quaternionic degree of symmetry} provides further evidence of the maximal symmetry rank conjecture.

We also point out that Galaz-García, Kerin and Radeschi proved that $\A(M)\leq \tfrac{2n}{3}$ if $M$ is a closed, simply connected rationally elliptic manifold, and they classify the manifolds for which the bound is attained up to equivariant rational homotopy equivalence (see \cite[Theorem A, Theorem B]{galaz2021torus}). They further conjectured that these manifolds admit an equivariant diffeomorphism classification (see \cite[Rigidity conjecture]{galaz2021torus}). \Cref{bound almost-effective quaternionic degree of symmetry} assumes the stronger symmetry of an almost-effective quaternionic torus action, but requires no rational ellipticity hypothesis and yields a classification up to diffeomorphism. Thus, \cref{bound almost-effective quaternionic degree of symmetry} also provides a quaternionic counterpart to the diffeomorphism conclusion of their rigidity conjecture.  

We also prove bounds and partial classification results for $\Q(M)$ and $\Qr(M)$, which are analogous to the results for $\A(M)$ and $\rank(M)$. 

\begin{thm}
\label{bound quaternionic degree of symmetry}
Let $M$ be a closed, connected smooth $n$-dimensional manifold. Then $\Q(M)\leq \frac{n}{3}$. Moreover, if $M$ is simply connected, $n=3k$ and $\Q(M)=\frac{n}{3}$, then $M$ is diffeomorphic to $Q^k$. 
\end{thm}

\begin{thm}
\label{bound quaternionic rank of symmetry}
Let $M$ be a closed, connected smooth $n$-dimensional manifold. Then $\Qr(M)\leq\tfrac{n}{3}$. Moreover:
\begin{itemize}
    \item[1.] If $n=3k$ and $\Qr(M)=\frac{n}{3}$, then $M$ is finitely covered by $Q^k$.
    \item[2.] If $n=3k+1$ and $\Qr(M)=\frac{n - 1}{3}$, then $M$ is finitely covered by $Q^k\times S^1$.
\end{itemize}
\end{thm}

A key feature of the almost-effective bound is that any action attaining it has cohomogeneity at most one. This permits the complete diffeomorphism classification in \cref{bound almost-effective quaternionic degree of symmetry}. By contrast, actions attaining the effective or almost-free dimension bounds can have cohomogeneity two, where classification is more difficult. Our classification results in these cases are therefore partial.

We also compute these invariants for some classes of manifolds.

\begin{thm}
\label{quaternionic degree of symmetry examples}
Let $M$ be a closed, connected smooth $n$-dimensional manifold. Then:
\begin{itemize}
    \item[1.] If $H_*(M,\Z/2)\cong H_*(S^n,\Z/2)$, then $\Q(M)\leq\tfrac{n+1}{4}$.
    \item[2.] If $\chi(M)\neq 0$, then $\Q(M)\leq\tfrac{n}{4}$.
    \item[3.] If $M\cong Q^a\times T^b$, then $\Q(M)=a$. 
\end{itemize}
\end{thm}

The first two bounds parallel the classical bounds $\A(S^n)\leq\tfrac{n+1}{2}$ and $\A(M)\leq\tfrac{n}{2}$ when $\chi(M)\neq0$. The third complements the computation $\A(Q^a\times T^b)=2a+b$ in \cite[Theorem C]{saito1984semisimple}.

\begin{cor}
\label{cor: quaternionic degree of symmetry examples}
Let $M$ be a closed, connected smooth $n$-dimensional manifold. Then:
\begin{itemize}
    \item[1.] If $M\cong S^n$, then $\Q(M)=\lfloor\tfrac{n+1}{4}\rfloor$.
    \item[2.] If $M$ is a quoric manifold of dimension $n=4k$, then $\Q(M)=\frac{n}{4}$.
    \item[3.] For any natural number $0\leq a\leq \lfloor\tfrac{n}{3}\rfloor$, there exists a closed, connected smooth $n$-dimensional manifold $M_a$ such that $\Q(M_a)=a$. 
\end{itemize}    
\end{cor}

The toral rank conjecture predicts that if $T^r$ acts effectively and almost-freely on a manifold $M$, then $\dim H^*(M,\QQ)\geq 2^r$ (see \cite{allday1993cohomological, felix2008algebraic} and references therein for background and known cases). In the spirit of the toral rank conjecture, we also obtain a cohomological lower bound when the quaternionic rank attains its maximal numerical value.

\begin{prop}
\label{prop: quaternioic rank maximal value}
Let $M$ be a closed, connected smooth $n$-dimensional manifold. If $\Qr(M)=\lfloor\tfrac{n}{3}\rfloor$, then
$$\dim H^*(M,\QQ)\geq 2^{\Qr(M)}.$$
\end{prop}

In a companion paper in preparation by the authors, we propose a quaternionic version of the toral rank conjecture and study it for quaternionic moment-angle manifolds.

\textbf{Organization of the paper.} The paper is structured as follows. The first section is the introduction. In the second section, we recall preliminaries on low-cohomogeneity actions and prove technical lemmas on connected closed subgroups of $Q^r$, which are essential to proving the results of this article. The third section is devoted to the quaternionic degree of symmetry, proving \cref{bound quaternionic degree of symmetry} and \cref{quaternionic degree of symmetry examples}. In the fourth section, we study the quaternionic rank of symmetry, and we prove \cref{bound quaternionic rank of symmetry}. Finally, in the fifth section, we prove \cref{bound almost-effective quaternionic degree of symmetry}.

\textbf{Funding.} The first author was supported by the research grant EXCELLENCE/0925-Hopf problem and monotone invariants, under the program “THALIA 2021-2027”, co-funded by the European Union through the Research and Innovation Foundation of Cyprus.

The second author was partially supported by the Spanish Ministry of Science, Innovation and Universities through the FPU fellowship (FPU22/03443).

The first and second authors were also partially supported by the grant PID2023-147642NB-I00 from the Spanish Ministry of Science, Innovation and Universities.

\textbf{AI disclosure.} All the ideas, arguments, results and proofs were obtained without the use of AI. We used ChatGPT for proofreading and improving the writing. The authors take full responsibility for the results.

\section{Preliminaries}

\subsection{Actions with low-cohomogeneity orbits}

We start by recalling some standard facts and terminology from the theory of compact transformation groups. For a detailed exposition, see \cite[Chapter I]{bredon1972introduction} or \cite[Chapter I]{lee1987manifolds}. 

Let $G$ be a compact Lie group acting on a closed, connected smooth manifold $M$. Given a point $x\in M$, the orbit $G(x)$ is a closed, connected smooth submanifold of $M$ which is $G$-equivariantly diffeomorphic to the coset space $G/G_x$. We say that two orbits $G/G_x$ and $G/G_y$ have the same \textit{type} if $G_x$ and $G_y$ are conjugate subgroups of $G$ or, equivalently, if there exists a $G$-equivariant diffeomorphism $G/G_x\longrightarrow G/G_y$. This defines an equivalence relation on the set of orbits of the $G$-action on $M$, and we denote its elements by $(G/H)$. We can define a partial order on the set of orbit types defined as follows: given two orbit types $(G/H)$ and $(G/K)$, we set $(G/H)\geq (G/K)$ if and only if $H$ is conjugate to a subgroup of $K$.  

Since $G$ is compact, there exists a maximal orbit type. An orbit belonging to this type is called \emph{principal}. The union of all principal orbits 
$$M_0=\{x\in M:\text{ $G(x)$ is a principal orbit}\} $$
is an open and dense subset of $M$ (see \cite[Chapter IV, Theorem 3.1]{bredon1972introduction}). A non-principal orbit is \emph{exceptional} if it has the same dimension as principal orbits, and \textit{singular} otherwise.

The \textit{cohomogeneity} of an action is the codimension of a principal orbit in $M$. For example, if the cohomogeneity of an action is zero, principal orbits have the same dimension as $M$ and are open submanifolds of $M$. Since principal orbits are also closed submanifolds and $M$ is connected, $M$ is diffeomorphic to a principal orbit. Then the action of $G$ on $M$ is transitive, hence $M\cong G/H$ is a homogeneous space. The following result on homogeneous spaces will be used throughout the paper. Although it is standard, we provide a proof for the convenience of the reader.

\begin{lem}
\label{lemma: homogeneous space effective}
Let $G/H$ be a compact homogeneous space with $G$ a compact, connected Lie group and $H$ a closed subgroup of $G$. The left transitive action of $G$ on $G/H$ is effective if and only if the normal core $\core_G(H)=\bigcap_{g \in G} gHg^{-1}$ is trivial. In particular, if $\Zc G\cap H$ is non-trivial, then the action of $G$ on $G/H$ is not effective.
\end{lem}

\begin{proof}
The kernel of the action of $G$ on $G/H$ is 
$$\Gamma = \{g \in G : gxH = xH \ \text{for all } x \in G \}.$$
It is clear that $\core_{G}(H) \subseteq \Gamma$. Conversely, if $g\in \Gamma$, then $x^{-1}gx\in H$ for all $x\in G$, hence $g\in \core_{G}(H)$ and $\Gamma \subseteq \core_{G}(H)$. Therefore, $\Gamma = \core_G(H)$, and the action is effective if and only if $\core_G(H)$ is trivial.

Let $h\in \Zc{G}\cap H$, then for all $x\in G$, we have $xhx^{-1}=h$, so $h\in xHx^{-1}$ for all $x\in G$. Consequently, $h\in\core_G(H)$ and $\Zc G\cap H\subseteq\core_G(H)$. Thus, if $\Zc G\cap H$ is non-trivial, then $\core_G(H)$ is non-trivial, and hence the action of $G$ on $G/H$ is not effective.
\end{proof}

We now recall some facts about cohomogeneity one actions on manifolds. First, we recall Mostert's result, which characterizes the structure of compact manifolds admitting cohomogeneity one actions.

\begin{thm}\cite{mostert1957compact},\cite[Chapter IV, \S 8]{bredon1972introduction}
\label{cohomogenity one structure}
    Let $M$ be a closed, connected, smooth manifold with a cohomogeneity one action of a compact, connected Lie group $G$, and let $G/H$ be a principal orbit. Then exactly one of the following holds:
	\begin{itemize}
		\item[1.] The orbit space is $M/G\cong S^1$. Then $M$ is $G$-equivariantly diffeomorphic to the total space of a smooth fibre bundle $\pi: M\longrightarrow S^1$ with fibre $G/H$ and structure group $N_G(H)/H$. The action has no non-principal orbits.
		\item[2.] The orbit space is $M/G\cong [-1,1]$. Then there exist subgroups $K_+$ and $K_-$ of $G$ containing $H$ and satisfying $K_\pm/H\cong S^{l_\pm}$ for some $l_{+}, l_{-} \geq 0$. The actions of $K_\pm$ on $S^{l_\pm}$ are linear and thus extend to linear actions of $K_\pm$ on $\mathbb{D}^{l_\pm+1}$. The non-principal orbits of the $G$-action on $M$ are $G / K_\pm$ and correspond to the preimages of $\pm 1 \in [-1,1]$, respectively. Moreover, $M$ is $G$-equivariantly diffeomorphic to the union of disk bundles over $G/K_\pm$ glued along the common boundary $G/H$, 
		$$M \cong \left( G\times_{K_+}\mathbb{D}^{l_++1} \right) \cup_{G/H} \left( G\times_{K_-}\mathbb{D}^{l_-+1} \right).$$
	\end{itemize}
\end{thm}

\begin{rmk}\label{cohomogeneity one diagram}
Conversely, let $G$ be a compact, connected Lie group and let $H$, $K_+$ and $K_-$ be closed subgroups whose inclusions form the commutative diagram
\[\begin{tikzcd}
    & G& \\
    K_+ \ar{ru}{\iota_+} & & K_-\ar{lu}[xshift=3ex,yshift=2.5ex]{\iota_-}\\
    & H\ar{lu}{\overline{\iota}_+}\ar{ru}[xshift=3ex,yshift=-3ex]{\overline{\iota}_-} &
\end{tikzcd}\]
Assume that $K_\pm/H\cong S^{l_\pm}$ for some $l_{+}, l_{-} \geq 0$, where these transitive actions are linear, and consider the extended action of $K_\pm$ on $\mathbb{D}^{l_\pm+1}$. The manifold 
$$\left( G\times_{K_+}\mathbb{D}^{l_++1} \right) \cup_{G/H} \left( G\times_{K_-}\mathbb{D}^{l_-+1} \right)$$
admits a cohomogeneity one action of $G$ with principal orbit $G/H$ and non-principal orbits $G/K_\pm$. Thus, the above group diagram characterizes cohomogeneity one actions on closed, connected manifolds whose quotient is an interval. We call this the \emph{group diagram} of the action and write it as $H\subseteq K_+, K_-\subseteq G$.
\end{rmk}

Different group diagrams can lead to manifolds which are $G$-equivariantly diffeomorphic.

\begin{prop}\cite[Proposition 1.2]{hoelscher2010classification}
\label{cohomogeneity one equivalent}
Let $M$ be a closed, connected manifold with a cohomogeneity one action of $G$ given by the group diagram $H \subseteq K_{+}, K_{-} \subseteq G$. Then any of the following operations on the group diagram will result in a $G$-equivariantly diffeomorphic manifold:
\begin{enumerate}
    \item[1.] switching $K_+$ and $K_-$,
    \item[2.] conjugating each group in the diagram by the same element of $G$, and
    \item[3.] replacing $K_-$ with $aK_-a^{-1}$ for some $a\in N(H)_{0}$.
\end{enumerate}
Conversely, the group diagrams for any two $G$-equivariantly diffeomorphic manifolds must be related by some combination of these three operations.
\end{prop}

\begin{rmk}
Let $M$ be a closed, connected manifold with a cohomogeneity one action of $G$ given by the group diagram $H \subseteq K_{+}, K_{-} \subseteq G$ with $K_\pm/H \cong S^{l_\pm}$ for some $l_{+}, l_{-} \geq 0$, and let $\phi: G\longrightarrow G$ be an automorphism. Then the group diagram $\phi(H) \subseteq \phi(K_+), \phi(K_-) \subseteq G$ satisfies $\phi(K_\pm)/\phi(H) \cong S^{l_\pm}$ and the resulting cohomogeneity one manifold is diffeomorphic to $M$, but it does not need to be $G$-equivariant diffeomorphic to $M$. 
\end{rmk}

We finish this section with a collection of lemmata that we use to prove the main results of this article.

\begin{lem}\cite[Theorem 3]{mostert1957compact}
\label{cohomogeneity one orientable}
Let $M$ be a closed, connected smooth manifold admitting a cohomogeneity one action of a compact, connected Lie group $G$ with a group diagram $H \subseteq K_+, K_- \subseteq G$ with $K_\pm/H \cong S^{l_\pm}$ for some $l_{+}, l_{-} \geq 0$. Assume that all the orbits are orientable. Then $M$ is orientable if and only if $l_\pm \geq 1$.   
\end{lem}

\begin{lem}\cite[Corollary 1.9]{hoelscher2010classification}
Let $M$ be a closed, connected smooth manifold with a cohomogeneity one action of a compact, connected Lie group $G$ with a group diagram $H \subseteq K_+, K_- \subseteq G$ with $K_\pm/H \cong S^{l_\pm}$ for some $l_{+}, l_{-} \geq 0$.
\begin{enumerate}
    \item[1.] If $l_+ > 1$ and $l_- \geq 1$, then $\pi_1(M) \cong \pi_1(G/K_-)$ and $H \cap K_+^0 = H^0$. In particular, if $M$ is simply connected, then $K_-$ is connected.
    \item[2.] If $l_-,l_+ > 1$, then $\pi_1(M) \cong \pi_1(G/H) \cong \pi_1(G/K_\pm)$. In particular, if $M$ is simply connected, then $H$, $K_-$, and $K_+$ are connected.
\end{enumerate}
\end{lem}

\begin{lem}\cite[1.21 Special type of actions]{hoelscher2010classification}
\label{cohomogeneity one product}
Let $M$ be a closed, connected, smooth manifold with a cohomogeneity one action of a compact, connected Lie group $G$ with a group diagram $H \subseteq K_+, K_- \subseteq G$ with $K_\pm/H \cong S^{l_\pm}$ for some $l_{+}, l_{-} \geq 0$. Let $G'/H'$ be a compact homogeneous space where $G'$ and $H'$ are compact Lie groups. Then $M\times G'/H'$ admits a cohomogeneity one action of $G\times G'$ with a group diagram $H\times H' \subseteq K_+\times H',K_-\times H' \subseteq G\times G'$ with $(K_\pm\times H')/(H\times H') \cong S^{l_\pm}$.

Conversely, every cohomogeneity one $(G\times G')$-manifold with this subgroup diagram, where all inclusions are the corresponding product inclusions, is $(G\times G')$-equivariantly diffeomorphic to $M\times G'/H'$.
\end{lem}

\subsection{Connected subgroups of \texorpdfstring{$Q^r$}{Qr}}

To study the three invariants presented in the introduction, we first need to understand the structure of connected Lie subgroups of $Q^r$. The following lemmata, proved using standard Lie group theory, generalize the results in \cite[Chapter 2]{Hopkinsonthesis}. For the standard facts about compact Lie groups and their representations used in this subsection, we refer to \cite{BrockerTomDieck}.

We define $[r]=\{1,\dots,r\}$. Given $I\subseteq [r]$, we define 
$$Q^I=\{(q_1,\dots,q_r)\in Q^r:q_j=e \text{ if } j\notin I \}.$$ 

\begin{lem}
\label{connected subgroups of Q^r}
Let $H$ be a closed, connected subgroup of $Q^r$, with an inclusion $\iota: H\longrightarrow Q^r$. Then, $H\cong Q^a\times T^b$. Moreover, there exist disjoint subsets $I_1,\dots, I_a,J\subseteq [r]$ satisfying $(\bigcup_{i=1}^a I_i)\cup J=[r]$ such that, up to conjugation by an element of $Q^r$, $\iota: H\longrightarrow Q^r$ satisfies:
\begin{itemize}
    \item[1.] For each $i\in [a]$, the restriction of $\iota$ to the $i$-th factor of $Q^a$ is the diagonal embedding $$\Delta_{I_i}:Q\longrightarrow Q^{I_i},$$ which maps $q\in Q$ to $\Delta_{I_i}(q)=(q_1,\dots q_r)$, where $q_j=q$ if $j\in I_i$ and $q_j=e$ if $j\notin I_i$.
    \item[2.] The inclusion $\iota$ maps $T^b$ to a subtorus of $Q^J$.
\end{itemize}
In particular, $a+b\leq r$.
\end{lem}

\begin{proof}
Since $H$ is a compact connected Lie group, we have that $H\cong (H_1\times\cdots \times H_a\times T^b)/\Gamma$, where $H_i$ are simply connected, simple, compact, connected Lie groups and $\Gamma$ is a finite subgroup of $\Zc(H_1\times\cdots \times H_a\times T^b)$. Let $\pi: H_1\times\cdots \times H_a\times T^b\longrightarrow H$ denote the quotient morphism. 

Our first claim is that $H_i\cong Q$ for all $i$. Indeed, for all $1\leq i\leq a$ and $1\leq j\leq r$, we consider the group morphisms $$f_{ij}:=p_j\circ \iota\circ\pi_{|H_i}:H_i\longrightarrow Q,$$ where $p_j:Q^r\longrightarrow Q$ is the projection to the $j$-th coordinate. Since $H_i$ is simple and $\ker f_{ij}$ is a normal subgroup of $H_i$, $\ker f_{ij}$ is either finite or $H_i$. Note that if $\ker f_{ij}=H_i$ for all $j$, then $\iota(\pi(H_i))$ is trivial, and since $\iota$ is injective, $\pi(H_i)$ is also trivial. This implies $H_i\leq \Gamma$, which is not possible since $\Gamma$ is finite. Consequently, there exists $j$ such that $\ker f_{ij}$ is finite. In particular, $\dim H_i\leq 3$. Since $H_i$ is simply connected, connected, and compact, we can conclude that $H_i\cong Q$. Thus, $f_{ij}: Q\longrightarrow Q$ is an endomorphism. By \cite[Proposition 2.0.3]{Hopkinsonthesis}, every non-trivial endomorphism $f_{ij}:Q\longrightarrow Q$ is an inner automorphism. Consequently, each $f_{ij}$ is either trivial or an inner automorphism. For each non-trivial $f_{ij}$, choose $q_{ij}\in Q$ such that $f_{ij}=C_{q_{ij}}$, and set $q_{ij}=1$ when $f_{ij}$ is trivial.

Similarly, the image of $\iota\circ\pi_{|T^b}:T^b\longrightarrow Q^r$ is contained in a maximal tori of $Q^r$. The injectivity of $\iota$ implies that $\ker \iota\circ\pi_{|T^b}=\ker \pi_{|T^b}$ and hence $\ker \iota\circ\pi_{|T^b}$ is finite.  

Let $I_i=\{j\in[r]: \text{ $\ker f_{ij}\neq H_i$}\}$ and let $J=[r]\setminus(\cup_{i=1}^aI_i)$. We have $\iota\circ\pi_{|H_i}=C_{x_i}\circ \Delta_{I_i}$, where $x_i=(q_{i1},\dots,q_{ir})\in Q^r$ satisfy that $q_{ij}$ is trivial if $j\notin I_i$ and $C_{x_i}:Q^r\longrightarrow Q^r$ denotes the conjugation by $x_i$. Our next claim is that the subsets $I_1,\dots, I_a, J\subseteq [r]$ are disjoint. To reach a contradiction, we assume that there exist $i$ and $i'$ such that $I_i\cap I_{i'}\neq \emptyset$. Let $k\in I_i\cap I_{i'}$. For all $q\in H_i$ and $q'\in H_{i'}$, we have the commuting relation  
\begin{align*}
\iota\circ\pi_{|H_i\times H_{i'}}(q,q')
&=
\bigl(\iota\circ\pi_{|H_i\times H_{i'}}(q,e)\bigr)\cdot
\bigl(\iota\circ\pi_{|H_i\times H_{i'}}(e,q')\bigr)\\
&=
\bigl(\iota\circ\pi_{|H_i\times H_{i'}}(e,q')\bigr)\cdot
\bigl(\iota\circ\pi_{|H_i\times H_{i'}}(q,e)\bigr).
\end{align*}

Therefore, in the $k$-th component we have 
$$(q_{ik}qq_{ik}^{-1})(q_{i'k}q'q_{i'k}^{-1})=(q_{i'k}q'q_{i'k}^{-1})(q_{ik}qq_{ik}^{-1})$$ for all $q\in H_i$ and $q'\in H_{i'}$. Since the conjugations $C_{q_{ik}}$ and $C_{q_{i'k}}$ are automorphisms of $Q$, the above equality would imply that $Q$ is commutative, reaching a contradiction. An analogous argument shows that $\iota\circ\pi(T^b)\leq Q^J$. Thus, the morphism $\iota\circ\pi: Q^a\times T^b\longrightarrow Q^r$ satisfies the conclusions of the lemma. 

The last step is to show that these properties descend to $\iota: H\longrightarrow Q^r$. Let $g=(q_1,\dots,q_a,t)\in\Gamma\leq Q^a\times T^b$. Then $\iota(\pi(g))$ is the identity of $Q^r$. Consequently, for all $1\leq i\leq a$, we have $C_{x_i}(\Delta_{I_i}(q_i))$ is trivial. All morphisms $C_{x_i}\circ \Delta_{I_i}$ are injective since they are a composition of an injective morphism and an automorphism. Thus, all $q_i$ are the identity element, and therefore $\Gamma$ injects into the $T^b$-factor of $Q^a\times T^b$. Recall that for any finite subgroup $F$ of $T^b$, $T^b/F\cong T^b$. Thus, $H=(Q^a\times T^b)/\Gamma =Q^a\times (T^b/\Gamma)\cong Q^a\times T^b$. Moreover, note that $\iota\circ\pi_{|Q^a}=\iota_{Q^a}$ which implies that the restriction of $\iota$ on the $i$-th factor is $C_{x_i}\circ\Delta_{I_i}$. Finally, since $\iota\circ \pi$ maps $T^b$ inside a maximal torus of $Q^J$, the inclusion $\iota$ also maps the $T^b$-factor of $H$ inside a maximal torus of $Q^J$.
\end{proof}

\begin{cor}
\label{subgroups containing maximal tori of Qr}
Let $H$ be a connected subgroup of $Q^r$ containing a maximal torus $T^r$. Then $H\cong Q^a\times T^b$ with $a+b=r$.
\end{cor}

\begin{proof}
By \cref{connected subgroups of Q^r}, $H\cong Q^a\times T^b$ with $a+b\leq r$. Thus, maximal tori of $H$ are of the form $T^{a+b}$. Since $T^r\leq H$, we have $r\leq a+b$ and hence $r=a+b$.
\end{proof}

\begin{lem}
\label{connected subgroups of Q^r almost-effective}
Let $H$ be a closed and connected subgroup of $Q^r$. The left transitive action of $Q^r$ on $Q^r/H$ is almost-effective if and only if $|I_i|\geq 2$ for all $i$, where $I_i$ are defined as in \cref{connected subgroups of Q^r}.  
\end{lem}

\begin{proof}
The action of $Q^r$ on $Q^r/H$ by left multiplication is almost-effective if and only if $\core_{Q^r}(H)$ is finite, or equivalently, the identity component of $\core_{Q^r}(H)$ is trivial. Thus, the action is not almost-effective if and only if there exists a non-trivial connected closed subgroup $H'$ of $H$ such that $H'$ is normal in $Q^r$. Thus, our first objective is to study closed connected normal subgroups of $Q^r$.

Let $G$ be a closed connected subgroup of $Q^r$ and let $\iota: G\longrightarrow Q^r$ denote the inclusion. By \cref{connected subgroups of Q^r}, $G\cong Q^a\times T^b$ for some natural numbers $a$ and $b$. Moreover, the description of $G$ and $\iota$ in \cref{connected subgroups of Q^r} implies that $G$ is normal in $Q^r$ if and only if each factor is normal in $Q^r$. No non-trivial connected abelian subgroup of $Q^r$ is normal, hence $b=0$. Thus, $G\cong Q^a$. By \cite[Proposition 2.1.4]{Hopkinsonthesis} (see also the observation at the end of page 17 of \cite{Hopkinsonthesis}), $G$ is normal if and only if $|I_i|=1$ for all $1\leq i\leq a$.

We will prove that the left transitive action of $Q^r$ on $Q^r/H$ is not almost-effective if and only if there exists an $i$ such that $|I_i|=1$. Indeed, if $|I_i|=1$ for some $i$ then the $i$-th component of $H=Q^a\times T^b$ is a normal subgroup of $Q^r$ and hence $Q\leq \core_{Q^r}(H)$, which implies that the action of $Q^r$ on $Q^r/H$ is not almost-effective. 

Conversely, if the action of $Q^r$ on $Q^r/H$ is not almost-effective, then $H=Q^a\times T^b$ contains a non-trivial connected subgroup $G$ which is normal in $Q^r$. Thus, $G=Q^{a'}$ and $G$ embeds into the $Q^a$ factor of $H$, thus $a'\leq a$. We can use \cref{connected subgroups of Q^r} to describe the inclusion $Q^{a'}\longrightarrow Q^a$, obtaining sets $I'_{i'}$ for all $1\leq i'\leq a'$. Similarly, using \cref{connected subgroups of Q^r}, we describe the inclusion $$Q^{a'}\longrightarrow Q^a\longrightarrow Q^r.$$ It is straightforward to show that for all $1\leq i'\leq a'$, the associated sets to this inclusion are 
$$I''_{i'}=\bigcup_{i\in I'_{i'}} I_i.$$

Thus, $|I''_{i'}|=\sum_{i\in I'_{i'}}|I_i|$. Since $G$ is normal in $Q^r$, we have $|I''_{i'}|=1$ for all $1\leq i'\leq a'$. As $|I_i|\geq 1$, we conclude that there exists $i\in[a]$ such that $|I_i|=1$.
\end{proof}

\begin{lem}
\label{irr rep of Q^axT^b}
The minimum dimension of a faithful real representation of $Q^a\times T^b$ is $4a+2b$. 
\end{lem}

\begin{proof}
Let $V$ be a real faithful representation of $Q^a\times T^b$ of minimal dimension. We can decompose it into irreducible representations 
$$V=\bigoplus_{i=1}^rV_i.$$
Each factor is a non-trivial real representation of $Q^a\times T^b$. Otherwise, we could remove the trivial factors, obtaining a real faithful representation of dimension strictly smaller than $\dim_\R V$, contradicting minimality.

Let $E=\Zc(Q^a\times T^b)[2]\cong (\Z/2)^{a+b}$ denote the $2$-torsion of the centre of $Q^a\times T^b$. Since each representation $V_i$ is irreducible, each element $x\in E$ acts on each $V_i$ as either $Id$ or $-Id$. Hence, for $1\leq i\leq r$ we can define characters $\epsilon_i:E\longrightarrow\Z/2$ satisfying that $\epsilon_i(x)=1$ if $x$ acts trivially on $V_i$ and $\epsilon_i(x)=-1$ if $x$ acts as $-Id$ on $V_i$. 

We consider the group morphism $\epsilon:E\longrightarrow (\Z/2)^r$ defined as $\epsilon=(\epsilon_1,\dots,\epsilon_r)$. Note that $\epsilon$ is injective. Indeed, if $x\in \ker \epsilon$, then $\epsilon_i(x)$ would act trivially on $V_i$ for all $i$ and hence $x$ would act trivially on $V$, contradicting the faithfulness hypothesis. Hence, $\epsilon$ is injective and $r\geq a+b$.

We consider the centre of $(\Z/2)^a\cong\Zc Q^a\leq E$. Let $k$ be the number of irreducible representations $V_i$ on which $\Zc Q^a$ acts non-trivially. We can reorder the factors of $V$ so that the irreducible representations on which $\Zc Q^a$ acts non-trivially are $V_1,\dots, V_k$. Since the restriction of $\epsilon_{|(\Z/2)^a}:(\Z/2)^a\longrightarrow (\Z/2)^r$ is also injective, there exist at least $a$ factors where an element of $\Zc Q^a$ (and hence $Q^a$) acts non-trivially. Thus, $k\geq a$. Given $x\in (\Z/2)^a$, there exists a subgroup of $Q^a$ isomorphic to $Q$ and with centre generated by $x$. Thus, each $V_i$ has a real faithful representation of $Q$ for $1\leq i \leq k$. Consequently, we have $\dim_\R V_{i}\geq 4$ for $1\leq i \leq k$. 

By minimality of $V$, each of the remaining factors, $V_{k+1},\dots V_{r}$, has a non-trivial real representation of $Q^a\times T^b$. Since $Q^a\times T^b$ is a compact, connected Lie group, every one-dimensional real representation is trivial. Consequently, $\dim_\R V_i\geq 2$ for all $k+1\leq i \leq r$.

By combining the two lower bounds of $\dim_\R V_i$, we obtain $$\dim_\R V\geq 4k+2(r-k)=2k+2r\geq 2a+2(a+b)=4a+2b.$$

On the other hand, the standard representation of $Q^a\times T^b$ on $V=\mathbb{H}^a\oplus\C^b$ is faithful, and it has real dimension $4a+2b$. In particular, since $V$ is a faithful real representation of minimal dimension, we have $\dim_\R V\leq 4a+2b$. Consequently, $\dim_\R V=4a+2b$, as desired. 
\end{proof}

\begin{rmk}\label{irr rep of Q^axT^b 2-torsion}
The same argument used in \cref{irr rep of Q^axT^b} shows that the minimum dimension of a real representation of $Q^a\times T^b$ which is faithful on the $\Zc(Q^a\times T^b)[2]$ is also $4a+2b$. 
\end{rmk}

\section{The quaternionic degree of symmetry}

\begin{proof}[Proof of \cref{bound quaternionic degree of symmetry}]
We proceed by induction on $n$. For $n=1$, we have $M \cong S^1$, so it admits no effective action of $Q$ by \cite[Corollary 3.1.17]{lee2010seifert}. Consequently, $\Q(M)=0$.

Assume now that $M$ is a closed, connected $n$-dimensional manifold and that the inequality holds for all closed, connected manifolds of dimension less than $n$. Assume $r=\Q(M)$, hence $Q^r$ acts effectively on $M$. Recall $M_0$ denotes the union of all principal orbits of the action. Then $Q^r$ acts effectively on $M_0$. This follows from $M_0$ being a dense invariant subset of $M$. 

This implies that $Q^r$ acts effectively on the principal orbits. Indeed, assume there exists a principal orbit $P=Q^r/H$ on which the action of $Q^r$ is not effective. By \cref{lemma: homogeneous space effective}, there exists a non-trivial element $g\in \core_{Q^r}(H)$. For any other principal orbit $P'=Q^r/H'$, the subgroups $H$ and $H'$ are conjugate, and hence $\core_{Q^r}(H)=\core_{Q^r}(H')$. Consequently, $g$ also acts trivially on $P'$. Thus, $g$ acts trivially on the union of all principal orbits $M_0$, which is a contradiction. Therefore, $Q^r$ acts effectively on a principal orbit $P=Q^r/H$.

By \cite[Chapter VI, Corollary 1.3]{bredon1972introduction}, $P$ is a closed submanifold of $M$. Assume that the codimension of $P$ is greater than or equal to $1$ in $M$. Then $\dim(P)\leq n-1$ and $\Q(M)\leq \Q(P)\leq\tfrac{n-1}{3}\leq \tfrac{n}{3}$ and we are done. 

Thus, it remains to prove the case where the codimension of $P$ is $0$. Let $H_0$ denote the connected component of the identity of $H$. Then $Q^r/H_0\longrightarrow Q^r/H$ is a regular $H/H_0$-covering, and since $Q^r$ is simply connected, we can lift any effective $Q^k$ action on $P$ to an effective action of $Q^k$ on $Q^r/H_0$ (see \cite[Chapter I, Theorem 9.1]{bredon1972introduction}).

We claim that $H_0$ must be trivial. Indeed, assume $H_0$ is non-trivial to obtain a contradiction. Since $H_0$ is compact and connected, it contains a subgroup isomorphic to $S^1$. Recall that the centre of $Q$ is isomorphic to $\Z/2$ and is generated by $-Id$. It is precisely the $2$-torsion subgroup of a maximal torus $S^1$ of $Q$. Hence, the centre $\Zc Q^r\cong(\Z/2)^r$ is the 2-torsion subgroup of any maximal torus $T^r$ of $Q^r$ (note that $\Zc Q^r$ is a subgroup of each maximal torus, since all maximal tori are conjugate and conjugation automorphisms fix the centre). The subgroup $S^1$ in $H_0$ is also contained in a maximal torus of $Q^r$, so the unique element $g \in S^1$ of order $2$ is in the centre of $\Zc Q^r$. This implies that $H_0\cap\Zc Q^r$ contains a non-trivial element. Therefore, by \cref{lemma: homogeneous space effective}, the action of $Q^r$ on $Q^r/H_0$ is not effective, which is a contradiction. Thus, $H_0$ has to be trivial. Since $H$ is a closed subgroup of $Q^{r}$, then $H$ is finite, so $n=\dim P = 3r$ because the codimension of $P$ is $0$. Therefore, $\Q(M)\leq\tfrac{n}{3}$, as claimed.

Assume now that $M$ is simply connected, $n=3k$, and $\Q(M)=k$. By the same arguments as above, $Q^k$ acts effectively on a principal orbit $P=Q^k/H$ of its action on $M$. Also, as before, $H_0$ needs to be trivial, and since $H$ is a closed subgroup of $Q^k$, $H$ is a finite subgroup of $Q^k$ with trivial intersection with its centre. Therefore, $\dim P=3k=\dim M$ and hence $P$ is open subset of $M$. Therefore, the connectedness of $M$ implies $P\cong M$. Applying the long exact sequence of homotopy to the covering $Q^k\longrightarrow Q^k/H$ yields $\pi_1(M)\cong H$. Finally, since $M$ is simply connected, $H$ is trivial, and hence $M\cong Q^k$.
\end{proof}

\begin{rmk}
Classifying manifolds of dimension $3k+1$ and $3k+2$ with $\Q(M)=k$ requires studying the possible orbits of the respective cohomogeneity one and cohomogeneity two actions, which is a hard problem in general. However, classification results exist in low dimensions, which we review below. 

First, if $M$ is a closed, connected smooth manifold and $\dim M$ is $1$ or $2$, then $\Q(M)=0$, so no closed connected manifold of dimension $1$ or $2$ admits an effective $Q$-action.

If $\dim M=4$ and $\Q(M)=1$, then $M$ admits a cohomogeneity one action of $Q$. Assume that the principal isotropy subgroup of this action is trivial. If $M/Q\cong S^1$, then $M$ is diffeomorphic to $S^3\times S^1$. If $M/Q\cong[-1,1]$, then $M$ is diffeomorphic to one of $S^4$, $\RP^4$, $\CP^2$, $\RP^4\#\RP^4$, $\RP^4\#\CP^2$, or $\CP^2\#(-\CP^2)$; see \cite[\S 2, cases 42--47]{parker1986}.

If $\dim(M)=5$ and $\Q(M)=1$, then the $Q$-action is of cohomogeneity $2$. These manifolds are classified in \cite{simas2012nonnegatively} with the extra assumption of simply-connectedness. In particular, $M$ is diffeomorphic to $S^5$, $S^2\times S^3$, the non-trivial $S^3$-bundle over $S^2$, or the Wu manifold $\SU(3)/\SO(3)$.

If $\dim(M)=7$ and $\Q(M)=2$, then $M$ admits a cohomogeneity one action of $Q^2$. Simply connected $7$-dimensional manifolds of cohomogeneity one are classified in \cite{hoelscher2010classification, hoelscher2009homology}.

For $\dim(M) = 8$ and $\Q(M) = 2$, we refer to \cite{mazaud2000Spin,mazaud20018dim}, which give a full equivariant classification of $8$-dimensional manifolds with effective $Q^2$-actions and principal orbits with trivial isotropy. 

No classification or general result is known in dimensions greater than or equal to $10$.
\end{rmk}


We prove now \cref{quaternionic degree of symmetry examples}, starting with the case of homology spheres.

\begin{prop}
\label{quaternionic degree of symmetry sphere}
Let $M$ be a closed, connected, orientable smooth $n$-dimensional manifold satisfying $H_*(M,\Z/2)\cong H_*(S^n,\Z/2)$. Then $\Q(M)\leq\tfrac{n+1}{4}$.    
\end{prop}

The key results to prove \cref{quaternionic degree of symmetry sphere} are the Borel formula and a lemma on invariant spin structures. We briefly recall some definitions of invariant spin structures. We refer to \cite{lawson2016spin,chichilnisky1972group,daura2022g} for general results on spin structures and invariant spin structures. 

Let $M$ be a closed, connected, oriented Riemannian manifold and let $FM$ denote the total space of the frame bundle of oriented orthonormal frames on $M$, which is a principal $\SO(n)$-bundle. Let $\lambda_n:\Spin(n)\longrightarrow\SO(n)$ denote the connected $2$-covering of $\SO(n)$. A spin structure is a tuple $(P,\Lambda)$ where $P\longrightarrow M$ is a principal $\Spin(n)$-bundle and $\Lambda:P\longrightarrow FM$ is a $2$-covering which is $\lambda_n$-equivariant. Recall that $M$ has a spin structure if and only if the second Stiefel-Whitney class vanishes, $w_2(M)=0$ (see \cite[Chapter II]{lawson2016spin}).

Assume $G$ is a compact Lie group acting effectively and by orientation-preserving isometries on $M$. Recall that we can make any metric $G$-invariant by averaging, so we do not lose any generality by assuming the action is by isometries. Then the $G$-action lifts to a free action on $FM$ such that $FM\longrightarrow M$ is $G$-equivariant. We say that a spin structure is $G$-invariant if there exists an action of $G$ on $P$ such that $\Lambda: P\longrightarrow FM$ is $G$-equivariant. The following two conditions characterize $G$-invariance of a spin structure $\Lambda:P\longrightarrow FM$ (see \cite[\S 5]{chichilnisky1972group}).
\begin{itemize}
    \item[1.] $G$ preserves the cohomology class of the $2$-covering associated to $\Lambda:P\longrightarrow FM$ in $H^1(FM,\Z/2)$. Then there exists a compact Lie group $\widetilde{G}$ fitting in the short exact sequence
        \[\begin{tikzcd}
            1\ar{r}{} & \Z/2 \ar{r}{} & \widetilde{G} \ar{r}{} & G \ar{r}{} & 1,
        \end{tikzcd}\] 
        where $\widetilde{G}$ acts on $P$ and $\Lambda$ is equivariant with respect to the morphism $\widetilde{G}\longrightarrow G$.
    \item[2.] There exists a section homomorphism $s:G\longrightarrow \widetilde{G}$ of the short exact sequence above.
\end{itemize}

Note that if $G$ is a connected Lie group, then the first condition is automatically satisfied. Moreover, if $G$ is simply connected, then $\widetilde{G}\cong G\times \Z/2$, so the second condition is also satisfied. In particular, if $Q^r$ acts on a spin manifold, then all its spin structures are $Q^r$-invariant. Note also that if $H\leq G$ and a spin structure is $G$-invariant, then it is also $H$-invariant.

\begin{lem}
\label{codimension spin structure}
Let $M$ be a closed, connected spin $n$-dimensional manifold. Assume that $\Z/2$ acts effectively on $M$ and that the following conditions hold:
\begin{itemize}
    \item[1.] there exists a spin structure $(P,\Lambda)$ which is $\Z/2$-invariant, and
    \item[2.] $M^{\Z/2}\neq\emptyset$. 
\end{itemize}
Then every connected component $F$ of $M^{\Z/2}$ satisfies $n-\dim F\equiv0\pmod4.$
\end{lem}

\begin{proof}
Let $f: M\longrightarrow M$ denote the involution which defines the $\Z/2$ action on $M$. By averaging a Riemannian metric, we can assume $f$ is an isometry. Let $df: FM\longrightarrow FM$ denote the lift of $f$ to the oriented frame bundle, and let $\overline{f}: P\longrightarrow P$ be the lift to the $\Z/2$-invariant spin structure, which satisfies $df\circ\Lambda=\Lambda\circ \overline{f}$ and $\overline{f}^2=Id_P$. 

Let $F$ be a connected component of $M^{\Z/2}$ and let $m=\dim F$ and $c=n-m$. Given $x\in F$, we can decompose $T_xM\cong T_xF\oplus\nu_x$ and we can choose an orthonormal basis $\{e_1,\dots e_m, e'_1,\dots,e_c'\}$ satisfying $d_xf(e_i)=e_i$ for all $1\leq i\leq m$ and $d_xf(e_j')=-e_j'$ for all $1\leq j \leq c$. In particular $d_xf=(Id_m,-Id_c)\in \SO(n)$. Since the spin structure is $\Z/2$-invariant, the set $\lambda_n^{-1}(d_xf)\in \Spin(n)$ cannot contain elements of order $4$. Otherwise, the section in condition 2 above could not exist. This condition happens if and only if $\tfrac{c}{2}$ is even (see \cite[Chapter I]{lawson2016spin}). Consequently, $c=n-\dim F =0\mod 4$.
\end{proof}

\begin{proof}[Proof of \cref{quaternionic degree of symmetry sphere}]
Assume that $Q^r$ acts effectively on $M$. Let $Z=\Zc Q^r$ denote the centre of $Q^r$ and let $z_j$ denote the involution generating the centre of the $j$-th factor of $Q^r$. Thus, $Z=\langle z_1,\cdots, z_r\rangle\cong (\Z/2)^r$. By Smith theory (see \cite[Chapter III, Theorem 5.1]{bredon1972introduction}), for every subgroup $H\leq Z$, the fixed-point set $M^H$ is either empty or a mod-$2$ homology sphere. In the latter case, it is either connected or consists of two points, so all its connected components have the same dimension. We therefore write $d(H)=\dim M^H$ if $M^H\neq\emptyset$ and $d(H)=-1$ otherwise. Since $Z$ is an elementary abelian $2$-group, we have $$n-d(Z)=\sum_{[Z:H]=2}(d(H)-d(Z))$$ by the Borel formula \cite[Chapter XIII, Theorem 2.3]{borel1960seminar}. 

Define $a_H=d(H)-d(Z)$. Our first claim is that if $a_H>0$, then $a_H\geq 4$. Since $[Z:H]=2$, there exists a number $j$ such that $z_j\notin H$. This follows from the fact that the elements $z_j$ generate $Z$. Note that $\langle z_j\rangle\cong \Z/2$, that $H$ and $z_j$ commute, and that $Z=\langle H,z_j\rangle$, which implies that $\langle z_j\rangle$ acts on the fixed point set $M^H$ and that $(M^H)^{\Z/2}=M^Z$. Since $a_H>0$, the action of $\langle z_j\rangle$ on $M^H$ is effective. Note that the $j$-th factor $Q_j$ of $Q^r$ and $H$ also commute, since $H$ is central in $Q^r$. This implies that $Q_j$ acts on $M^H$, and this action is effective since the action of its centre $\langle z_j\rangle$ is effective. Since $Q_j$ is connected, it cannot act effectively on two points. Thus, $M^H$ is connected. In particular, $d(H)\geq 3$ by \cref{bound quaternionic degree of symmetry}. By \cite[Chapter III, Theorem 5.1]{bredon1972introduction}, $H_*(M^H,\Z/2)\cong H_*(S^{d(H)},\Z/2)$. Since $d(H)\geq 3$, we have $H^1(M^H,\Z/2)=H^2(M^H,\Z/2)=0$ and hence the first and second Stiefel-Whitney classes $w_1(M^H)=w_2(M^H)=0$. Therefore, $M^H$ is spin. Since $Q_j$ is simply connected, the spin structures of $M^H$ are $Q_j$-invariant. Consequently, the spin structures of $M^H$ are also $\langle z_j\rangle$-invariant. If $M^Z\neq \emptyset$, then we can apply \cref{codimension spin structure} to conclude that $a_H\geq 4$. If $M^Z=\emptyset$, then $d(Z)=-1$ and therefore $a_H=d(H)-d(Z)\geq 3+1=4$. Thus, in both cases, $a_H\geq 4$. 

Our next goal is to show that there are at least $r$ subgroups $H\leq Z$ of index $[Z:H]=2$ with $a_H>0$. Indeed, let $\mathcal{H}=\{H\leq Z: \text{$[Z:H]=2$ and $a_H>0$}\}$ and let $A=\bigcap_{H\in \mathcal{H}} H$. We claim that $A$ is trivial. We consider the action of $Z/A$ on $M^A$. Subgroups of index $2$ in $Z/A$ correspond to subgroups of index $2$ of $Z$ containing $A$. Moreover, any subgroup $H\leq Z$ containing $A$ satisfies $(M^A)^{H/A}=M^H$. Consequently, $\dim (M^A)^{H/A}=d(H) $ and 
$$d(A)-d(Z)=\sum_{[Z/A:H/A]=2}(\dim (M^A)^{H/A}-d(Z))=\sum_{[Z:H]=2, A\leq H}a_H.$$ 
Since the sum is over all $a_H$ satisfying $a_H>0$, we have $d(A)-d(Z)=n-d(Z)$ and therefore $n=d(A)$. Consequently, $M^A$ has a connected component that is an $n$-dimensional submanifold, so it is an open submanifold of $M$. Since $M^A$ is also a closed submanifold and $M$ is connected, we can conclude that $M^A=M$, which implies that $A$ is trivial since the action of $Z$ on $M$ is effective. We can look at $Z$ as a $\Z/2$-vector space and each subgroup $H$ as a hyperplane of $Z$. Recall that given $s$ hyperplanes $H_1,\dots,H_s$, we have $\dim (H_1\cap\cdots\cap H_s)\geq r- s$. Since $A$ is trivial, we can conclude that $|\mathcal{H}|\geq r$, as claimed.

Finally, we can apply Borel formula again to obtain $$n+1\geq n-d(Z)=\sum_{[Z:H]=2}a_H\geq 4r.$$ Consequently, $r\leq\tfrac{n+1}{4}$ and hence $\Q(M)\leq\tfrac{n+1}{4}$.
\end{proof}


\begin{prop}
\label{quaternionic degree of symmetry quoric manifolds}
Let $M$ be a closed, connected smooth $n$-dimensional manifold. Assume that $\chi(M)\neq 0$. Then $\Q(M)\leq \tfrac{n}{4}$.
\end{prop}

\begin{proof}
Assume that $Q^r$ acts effectively on $M$ and let $T^r$ be a maximal torus of $Q^r$. We have $\chi(M^{T^r})=\chi(M)\neq 0$ and hence $M^{T^r}\neq \emptyset$, see \cite[Theorem 7.33]{felix2008algebraic}. Choose $x\in M^{T^r}$ and consider its orbit $P\cong Q^r/Q^r_x$. The connected component of $Q^r_x$ is $Q^a\times T^{r-a}$ by \cref{subgroups containing maximal tori of Qr} and thus $$\dim P=3r-3a-(r-a)=2(r-a).$$

We have the $Q_x^r$-invariant decomposition $T_xM\cong T_xP\oplus \nu_x$. We consider $(\Z/2)^r\cong\Zc Q^r\leq T^r\leq Q_x^r$, which acts faithfully on $T_xM$ because the action of $Q^r$ on $M$ is effective. Moreover, since $\Zc Q^r$ commutes with $Q_x^r$, the action of $\Zc Q^r$ on $P$ is trivial (see the proof of \cref{lemma: homogeneous space effective}), and hence trivial on $T_xP$. Since $T_xM\cong T_xP\oplus \nu_x$ is also $\Zc Q^r$-invariant, we obtain that $\Zc Q^r$ acts faithfully on $\nu_x$. Hence, the 2-torsion of the centre of the connected component $Q_x^r$ also acts faithfully on $\nu_x$. \Cref{irr rep of Q^axT^b} and \cref{irr rep of Q^axT^b 2-torsion} imply that $\dim\nu_x\geq 4a+2(r-a)$. Thus, we have
$$ n=\dim T_xM=\dim T_xP+\dim \nu_x=2(r-a)+\dim\nu_x\geq 2(r-a)+4a+2(r-a)=4r.$$ Consequently, $\Q(M)\leq \tfrac{n}{4}$. 
\end{proof}


\begin{lem}
\label{quaternionic degree of symmetry QxT}
For any numbers $a,b\geq 0$, we have $\Q(Q^a\times T^b)=a$. 
\end{lem}

The main ingredient in the proof of \cref{quaternionic degree of symmetry QxT} is the following result of Burghelea and Schultz:

\begin{prop}\cite[Proposition 2.b)]{burghelea1975semisimple}
\label{obstruction cup product}
Let $G$ be a compact Lie group acting on a closed, connected $n$-dimensional manifold $M$. Let $k$ be the dimension of a principal orbit. Assume that for any point $x\in M$ the morphism $$i_x^*:H^k(M,\R)\longrightarrow H^k(G/G_x,\R) $$ is trivial. Then for any cohomology classes $\gamma_1,\dots\gamma_{n-k}\in H^1(M,\R)$ and $\gamma\in H^k(M,\R)$ we have
$$\gamma_1\smile\cdots\smile\gamma_{n-k}\smile\gamma=0.$$
\end{prop}

\begin{proof}[Proof of \cref{quaternionic degree of symmetry QxT}]
Assume that $Q^r$ acts effectively on $Q^a\times T^b$ and that $r>a$. Since the action is effective, the isotropy subgroups of principal orbits are finite, and principal orbits have dimension $3r$. 

The cohomology ring of $H^*(M,\R)$ is generated by $\alpha_1,\dots,\alpha_a$ of degree $3$ and $\beta_1\dots,\beta_b$ of degree 1, corresponding to the generator of the cohomology of each factor $Q$ and $T$, respectively. In particular, 
$$\alpha_1\smile\cdots\alpha_a\smile \beta_1\smile \cdots\smile \beta_b\neq 0.$$

Since $r>a$, we have $3a+b-3r\leq b$. In particular, we can choose $3a+b-3r$ cohomology classes of degree $1$, $\gamma_i=\beta_i$ for all $1\leq i\leq 3a+b-3r$, and a class $\gamma=\beta_{3a+b-3r+1}\smile\cdots\smile \beta_b\smile \alpha_1\smile\cdots\smile \alpha_a$ satisfying $\gamma_1\smile\cdots\smile\gamma_{n-k}\smile\gamma\neq0$. By \cref{obstruction cup product}, there exists $x\in M$ such that $i_x^*:H^{3r}(M,\R)\longrightarrow H^{3r}(Q^r/Q^r_x,\R)$ is not trivial. Since $\dim Q^r/Q^r_x= 3r-\dim Q_x^r\leq 3r$, we necessarily have that $Q_x^r$ is finite. Therefore $H^{3r}(Q^r/Q^r_x,\R)\cong H^{3r}(Q^r,\R)^{Q_x^r}\cong H^{3r}(Q^r,\R)$, since the action of $Q_x^r$ on $H^*(Q^r,\R)$ is the restriction of the action of $Q^r$ on $H^*(Q_r,\R)$, which is trivial since $Q^r$ is connected. 

Note that $i_x^*(\beta_i)=0$ for all $1\leq i\leq b$ since $H^*(Q^r,\R)$ does not have elements of degree 1. Furthermore, every element $\gamma\in H^{3r}(M,\R)$ is the linear combination of elements of the form $\alpha_{i_1}\smile \cdots\smile\alpha_{i_p}\smile\beta_{j_1}\smile\cdots\smile\beta_{j_q}$ where $3p+q=3r$. Since $r>a$, we necessarily have that $q>0$. This implies that 
\begin{align*}
&i_x^*(\alpha_{i_1}\smile\cdots\smile\alpha_{i_p}
\smile\beta_{j_1}\smile\cdots\smile\beta_{j_q})\\
&\quad= i_x^*\alpha_{i_1}\smile\cdots
\smile i_x^*\alpha_{i_p}\smile i_x^*\beta_{j_1}\smile\cdots
\smile i_x^*\beta_{j_q}\\
&\quad=0.
\end{align*}

Consequently, $i_x^*:H^{3r}(M,\R)\longrightarrow H^{3r}(Q^r/Q^r_x,\R)$ is trivial, reaching a contradiction. We have proved that $\Q(Q^a\times T^b)\leq a$. Note that $Q^a$ acts freely on $Q^a\times T^b$ by left multiplication on the $Q^a$ factor and trivially on the $T^b$ factor. Hence \(a\leq\Q(Q^a\times T^b)\). Combining this with the preceding upper bound, we obtain \(\Q(Q^a\times T^b)=a\).
\end{proof}

We end this section with the proof of \cref{cor: quaternionic degree of symmetry examples}

\begin{proof}[Proof of \cref{cor: quaternionic degree of symmetry examples}]
First, assume that $M\cong S^n$. Then $\SO(n+1)$ acts effectively on $S^n$. Let $k=\lfloor \tfrac{n+1}{4}\rfloor$. Since $\SO(4k)\leq\SO(n+1)$ and $Q^k\leq\SO(4k)$, the group $Q^k$ acts effectively on $S^n$. Consequently, $k\leq\Q(S^n)$. On the other hand, \cref{quaternionic degree of symmetry examples}.1 implies $\Q(S^n)\leq \left\lfloor\frac{n+1}{4}\right\rfloor=k$. Therefore, $\Q(S^n)=k$.

We assume now that $M$ is a quoric manifold of dimension $n=4k$. By \cite[Theorem 4.2.12]{Hopkinsonthesis}, we have $\chi(M)\neq 0$ and thus $\Q(M)\leq\frac{n}{4}=k$ by \cref{quaternionic degree of symmetry examples}.2. On the other hand, $M$ admits an effective $Q^k$ action by definition, and hence $k \leq \Q(M)$, which implies that $\Q(M)=k$.

We now prove the third claim. Given a natural number $0\leq a\leq \lfloor\tfrac{n}{3}\rfloor$, we define $M_a\cong Q^a\times T^{n-3a}$. Then $M_a$ is a closed, connected smooth $n$-dimensional manifold satisfying $\Q(M_a)=a$ by \cref{quaternionic degree of symmetry examples}.3.
\end{proof}

\section{The quaternionic rank}

\begin{proof}[Proof of \cref{bound quaternionic rank of symmetry}]

The argument to obtain the bound follows the same strategy as in the proof of \cref{bound quaternionic degree of symmetry}. We proceed by induction on $n$. For $n=1$, $M\cong S^1$, and therefore admits no almost-free action of $Q$ by \cite[Corollary 3.1.17]{lee2010seifert}. Consequently, $\Qr(M)=0$.

Assume now that $M$ is a closed, connected $n$-dimensional manifold and that the inequality holds for all closed, connected manifolds of dimension less than $n$. Assume $r=\Qr(M)$. Thus, $M$ admits an effective $Q^r$-action. As in the proof of \cref{bound quaternionic degree of symmetry}, $Q^r$ acts effectively on a principal orbit $P=Q^r/H$. Since the action of $Q^r$ on $M$ is almost-free, $H$ is finite, so $Q^r$ acts almost-freely on $P$. If $P$ is a submanifold of codimension greater than or equal to $1$ in $M$, then $\dim(P)\leq n-1$. Thus, $\Qr(M)\leq \Qr(P)\leq\tfrac{n-1}{3}\leq \tfrac{n}{3}$ and we are done. 

If the codimension of $P$ is $0$, then $n = \dim P = 3r$, hence $\Qr(M)\leq \tfrac{n}{3}$, as claimed.

Assume now that $n=3k$ and $\Qr(M)=k$. By the same reasoning as before, $Q^{k}$ acts effectively and almost-freely on a principal orbit $P = Q^{k}/\Gamma$, where $\Gamma$ is finite. Then $P$ has codimension $0$, so it is open and closed in $M$. Since $M$ is connected, we obtain that $M \cong P \cong Q^{k}/\Gamma$. In particular, $M$ is finitely covered by $Q^k$. 

We assume now that $n=3k+1$ and $\Qr(M)=k$. Then $M$ admits a $Q^{k}$-action with principal orbits of codimension $1$. Assume that $M/Q^k\cong [-1,1]$. All isotropy subgroups are finite since the action is almost-free. Hence the orbits are of the form $Q^k/\Gamma$, where $\Gamma$ is a finite subgroup of $Q^k$. In particular, all orbits have dimension $3k$ and are orientable. By \cref{cohomogeneity one orientable}, $M$ is not orientable. Consider the orientable $2$-cover $p':M'\longrightarrow M$. Since $Q^k$ is simply connected, we can lift the action of $Q^k$ on $M$ to an action of $Q^k$ on $M'$ such that the covering $p'$ is $Q^k$-equivariant (see \cite[Chapter I, Theorem 9.1]{bredon1972introduction}). Thus, the $Q^k$-action on $M'$ is also almost-free and of cohomogeneity one. Since $M'$ is orientable, $M'/Q^k$ cannot be an interval. Therefore $M'/Q^k\cong S^1$ and we have a fibre bundle
$$Q^k/\Gamma \longrightarrow M' \longrightarrow S^1$$ 
with structural group $N_{Q^k}(\Gamma)/\Gamma$.

Using the long exact sequence of homotopy and the fact that $Q^k$ is simply connected, we obtain a short exact sequence
$$1\longrightarrow \Gamma\longrightarrow \pi_1(M')\longrightarrow \Z\longrightarrow 1.$$ 

In particular, there exists $\phi\in\Aut(\Gamma)$ such that $\pi_1(M')\cong\Gamma\rtimes_\phi\Z$. Since $\Gamma$ is finite, $\Aut(\Gamma)$ is also finite, so there exists an integer $a$ such that $\phi^a=Id_\Gamma$. Thus, $\Gamma\times\Z\cong \Gamma\rtimes_{\phi^a}\Z$ is a finite index subgroup of $\pi_1(M')$. We can pass to a finite cover $M''$ such that $\pi_1(M'')\cong \Gamma\times \Z$. Finally, there exists a finite cover $M'''$ of $M''$ such that $\pi_1(M''')\cong \Z$. Thus, we obtain a finite cover $q: M'''\longrightarrow M$. We can use \cite[Chapter I, Theorem 9.1]{bredon1972introduction} to lift the action of $Q^k$ on $M$ to an action of $Q^k$ on $M'''$ such that $q$ is $Q^k$-equivariant. Consequently, the action of $Q^k$ on $M'''$ is almost-free and of cohomogeneity one. Therefore, $M'''/Q^k\cong S^1$ and we have a fibre bundle
$$Q^k/\Gamma'''\longrightarrow M'''\longrightarrow S^1$$ 
with structural group $N_{Q^k}(\Gamma''')/\Gamma'''$. The long exact sequence of homotopy groups shows that $\Gamma'''$ is a subgroup of $\pi_1(M''')\cong \Z$. Since $\Gamma'''$ is finite, it is trivial, and hence $M'''\longrightarrow S^1$ is a principal $Q^k$-bundle. These principal bundles are classified by $[S^1, BQ^k]$, which is trivial since $Q^k$ is connected. Consequently, $M'''\cong Q^k\times S^1$.

Assume now that $M/Q^k\cong S^1$. Then $M$ is the total space of a fibre bundle $Q^k/\Gamma\longrightarrow M\longrightarrow S^1$ with structure group $N_{Q^k}(\Gamma)/\Gamma$. Applying the same finite-cover argument as above, we conclude that $M$ is finitely covered by $Q^k\times S^1$.
\end{proof}

\begin{rmk}    
Effective almost-free actions of $\SU(2)$ on $5$-dimensional manifolds were also  classified in \cite{simas2012nonnegatively}. In particular, if $M$ is a simply-connected $5$ dimensional manifold with an effective almost-free action of $\SU(2)$ then $M$ is diffeomorphic to $S^3\times S^2$ or the non-trivial $S^3$-bundle over $S^2$ (see \cite[Theorem C]{simas2012nonnegatively}).
\end{rmk}

For example, we can use a standard spectral sequence argument to compute the quaternionic rank for spheres.

\begin{prop}
\label{prop:sphere-rank}
$\Qr(S^n)=1$ if and only if $n=4k-1$ for a positive integer $k$, and $\Qr(S^n)=0$ otherwise.    
\end{prop}

\begin{proof}
Assume $Q^{r}$ acts effectively and almost-freely on $S^n$ with $r\geq 1$. Consider the Serre spectral sequence with rational coefficients of the Borel fibration $S^n\longrightarrow S^n_{Q^r}\longrightarrow BQ^r$,
$$ E_2^{i,j}=H^i(BQ^r,H^j(S^n,\QQ))\implies H^{i+j}(S^n_{Q^r},\QQ).$$
The second page is concentrated in the rows $0$ and $n$. Thus all $d_k:E^{i,j}_k\longrightarrow E^{i+k,j-k+1}_k$ are trivial if $k\neq n+1$. In particular, $E_2\cong E_{n+1}$ and $E_{n+2}\cong E_\infty$. Moreover, $H^*(BQ^r,\QQ)\cong \QQ[x_1,\dots,x_r]$, where $\deg(x_i)=4$ for all $i$, therefore $E^{i,j}_k=0$ for all $i\neq 0,4a$ with $a$ a positive integer.

The cohomology of the Borel construction $H^*(S^n_{Q^r},\QQ)$ is finite dimensional by \cite[Theorem 7.7]{felix2008algebraic}. This implies that $d_{n+1}:E^{4a,n}\longrightarrow E^{4a+n+1,0}$ needs to be an isomorphism for $a$ large enough. In particular, $E^{4a+n+1,0}\neq0 $ and therefore $4a+n+1=4b$ for some integer $b$ with $b > a$. This implies that $n=4(b-a)-1$.

Since $H^*(BQ^r,\QQ)\cong \QQ[x_1,\dots,x_r]$, we have $\dim E^{4a,0}_2= \dim E_2^{4a,n}=\binom{r-1 + a}{r-1}$. Thus, we have $d_{n+1}:E^{4a,n} \cong \QQ^{\binom{r-1+a}{r-1}} \longrightarrow E^{4a+n+1,0}\cong \QQ^{\binom{r-1+b}{r-1}}$. Since $b>a$, we have $\binom{r-1+a}{r-1}=\binom{r-1+b}{r-1}$ if and only if $r-1=0 $.  Thus, $d_{n+1}$ is an isomorphism only if $r=1$. 

We have proved that $\Qr(S^n)=0$, unless $n=4k-1$ for some positive integer $k$. In the latter case, $\Qr(S^n) \leq 1$. Finally, note that $S^3$ acts freely on $S^{4k-1}$ and $S^{4k-1}/S^3\cong \mathbb{H}P^{k-1}$. Therefore, $\Qr(S^{4k-1})=1$ for all positive integers $k$.
\end{proof}

\begin{rmk}
Note that if $\chi(M)\neq 0$, then $\Qr(M)=0$. Indeed, if $Q^r$ acts almost-freely on $M$, then $T^r$ also acts almost-freely on $M$, and therefore $\chi(M)=0$ by \cite[Theorem 7.33]{felix2008algebraic}.

If $M\cong Q^a\times T^b$, then $\Qr(M)=a$. Indeed, $\Qr(M)\leq \Q(M)=a$. Moreover, the action of $Q^a$ by left multiplication on the $Q^a$ factor and trivial on the $T^b$ factor is free. Hence $a\leq \Qr(M)$, and thus $\Qr(M)=a$.
\end{rmk}

We finish this section, bounding the cohomology of $M$ when $\Qr(M)=\lfloor\tfrac{n}{3}\rfloor$, proving \cref{prop: quaternioic rank maximal value},

\begin{proof}[Proof of \cref{prop: quaternioic rank maximal value}]
Let $r=\Qr(M)=\lfloor\tfrac{n}{3}\rfloor$ and consider the principal $Q^r$-bundle $$Q^r\longrightarrow M\times EQ^r\longrightarrow M_{Q^r}.$$ 
Since the action of $Q^r$ is almost-free, the natural map $M_{Q^r}\longrightarrow M/Q^r$ induces an isomorphism $H^*(M_{Q^r},\QQ)\cong H^*(M/Q^r,\QQ)$, see \cite[Theorem 7.6(2)]{felix2008algebraic}. Moreover, $M/Q^r$ is compact and admits the structure of an orbifold of dimension $n-3r\leq 2$, see \cite[\S 2.1]{galaz2021torus}. By \cite[Example 2.11, Proposition 2.12, \S 2.1 pp. 32-34]{adem2007orbifolds}, $H^i(M_{Q^r},\QQ)=0$ for all $i>2$ (see \cite[\S 2.1]{adem2007orbifolds} for an introduction to cohomology of orbifolds). 

Since $Q^r$ is a connected Lie group and the fibre bundle is principal, the action of $\pi_1(M_{Q^r})$ on the fibres is by translations, and hence $\pi_1(M_{Q^r})$ acts trivially on the cohomology $H^j(Q^r,\QQ)$ for all $j\geq 0$. Thus, the system of local coefficients of the Serre spectral sequence is trivial and hence $$E^{i,j}_2\cong H^i(M_{Q^r},\QQ)\otimes H^j(Q^r,\QQ)\implies H^{i+j}(M,\QQ).$$

Since $H^*(Q^r,\QQ)\cong \bigwedge(x_1,\dots,x_r)$ satisfying $\deg(x_k)=3$ for all $1\leq k\leq r$, $H^*(Q^r,\QQ)$ only has non-trivial cohomology in degrees divisible by $3$. This implies that all the differentials $d_2:E^{i,j}_2\longrightarrow E_2^{i+2,j-1}$ are trivial. For any $a\geq 3$, all the differentials $d_a:E^{i,j}_a\longrightarrow E^{i+a,j-a+1}_a$ are trivial, since $H^i(M_{Q^r},\QQ)=0$ for all $i>2$. Consequently, the Serre spectral sequence collapses at the second page, and therefore we have an isomorphism of $\QQ$-vector spaces $$H^*(M,\QQ)\cong H^*(M\times EQ^r,\QQ)\cong H^*(Q^r,\QQ)\otimes H^*(M_{Q^r},\QQ).$$

Since $\dim H^*(Q^r,\QQ)=2^r$ and $\dim H^*(M_{Q^r},\QQ)\geq 1$, we obtain
$$\dim H^*(M,\QQ)=\dim H^*(Q^r,\QQ)\dim H^*(M_{Q^r},\QQ)\geq 2^r,$$ as desired.
\end{proof}

\section{Almost-effective quaternionic degree of symmetry}

Recall that given an almost-effective group action of a compact, connected, Lie group $G$ on a manifold $M$
$$\phi: G\longrightarrow \Diff(M),$$
the finite group $\ker\phi$ is called the \emph{kernel} of the action.

\begin{proof}[Proof of \cref{bound almost-effective quaternionic degree of symmetry}]  

The proof follows the same strategy as the other two cases. However, the classification part is substantially longer. Thus, we divide the proof into seven parts.

\paragraph*{\textbf{Part 1: Bounding $\Q_{ae}$.}} As in the previous results, we proceed by induction on the dimension of $M$. If $n=1$, then $M\cong S^1$. Assume that $Q^r$ acts almost-effectively on $S^1$ with kernel a finite group $\Gamma$. Then $Q^r/\Gamma$ acts effectively on $S^1$. Since the only compact connected Lie group which acts effectively on $S^1$ is $S^1$ itself, we can conclude that $r=0$ and hence $\Q_{ae}(M)=0$.

Assume $M$ is a closed, connected $n$-dimensional manifold and that the inequality holds for all closed, connected manifolds of dimension less than $n$. Let $r=\Q_{ae}(M)$, and hence there exists an almost-effective action of $Q^r$ on $M$. Let $\Gamma$ be the kernel and let $G=Q^r/\Gamma$. Thus, $G$ acts effectively on $M$. Using the same arguments as in \cref{bound quaternionic degree of symmetry}, we can conclude that $G$ acts effectively on a principal orbit $P\cong G/H$. Using the quotient morphism $\pi:Q^r\longrightarrow G$ we can construct an almost-effective action of $Q^r$ on $G/H$ given by $q(gH)=(\pi(q)g)H$ for all $q\in Q^r$ and $gH\in G/H$. If the codimension of $P$ is greater than or equal to $1$, then $\Q_{ae}(M)\leq\Q_{ae}(P)\leq \tfrac{2(n-1)}{3}\leq \tfrac{2n}{3}$. 

We now prove the case where the codimension of $P$ is $0$. Then $P$ is an open and closed submanifold of $M$. Thus, $M$ is connected, so $M\cong G/H$. Using the projection $\pi: Q^r\longrightarrow G$, we can construct a transitive almost-effective action of $Q^r$ on $P$. Therefore, $P\cong Q^r/H'$, where $H'=\pi^{-1}(H)$. Let $H_0'$ denote the connected component of the identity of $H'$. Since $\core_{Q^r}(H')$ is finite and $\core_{Q^r}(H_0')\leq \core_{Q^r}(H')$, the action of $Q^r$ on $Q^r/H_0'$ is also transitive and almost-effective. Thus, we have $$n=\dim P=\dim Q^r-\dim H'_0 =3r-\dim H_0'.$$ 

From now on we use the same notation as in \cref{connected subgroups of Q^r}. We have $H_0\cong Q^a\times T^b$ for some natural numbers $a$ and $b$. Hence $\dim H_0'=3a+b$, and therefore $3r-n=3a+b$. We also have $$r=\sum_{i=1}^a|I_i|+|J|.$$ Since $|I_i|\geq 2$ for all $i$ and $|J|\geq b$, we have $$r\geq 2a+b.$$
Consequently, $$6r-2n+b=6a+3b\leq 3r,$$
which implies that $$r\leq\frac{2n-b}{3}\leq \frac{2n}{3},$$
since $b\geq 0$. Thus, $\Q_{ae}(M)=r\leq \tfrac{2n}{3}$, as desired.

\paragraph*{\textbf{Part 2: Rigidity for $n=3k$.}} We start now proving the rigidity result in each of the three cases. We start with the case where $n=3k$, hence $\Q_{ae}(M)=2k$. Let $Q^{2k}$ act almost-effectively on $M$. Then $Q^{2k}$ also acts almost-effectively on a principal orbit $P$, and since $\tfrac{3}{2}\Q_{ae}(M)\leq n$, we can conclude that $\dim P=n$ and hence $M\cong Q^{2k}/H$. Moreover, since $Q$ is simply connected, the long exact sequence of homotopy yields $\pi_1(M)\cong\pi_0(H)$. Thus, the simply-connectedness of $M$ implies that $H$ is connected. By \cref{connected subgroups of Q^r} $H\cong Q^a\times T^b$. We have
$$3k=\dim M= \dim Q^{2k}-\dim H=6k-3a-b$$
and therefore $$3k=3a+b.$$

We have the inequality $2k\geq 2a+b$. Combining both expressions, we obtain $a\geq k$. Since the action is almost-effective, the set $I_i$ satisfy $|I_i|\geq 2$ for all $1\leq i\leq a$ by \cref{connected subgroups of Q^r almost-effective}, which implies that $a\leq k$. Thus, $a=k$ and hence $b=0$. We can conclude that $H\cong Q^k$ and that $|I_i|=2$ for all $i$. By \cref{connected subgroups of Q^r}, we can modify the action of $Q^{2k}$ by precomposing it with an inner automorphism and a permutation of the factors of $Q^{2k}$ so the sets $I_i$ satisfy $I_{i}=\{2i-1,2i\}$ for all $1\leq i\leq a$ and the inclusion of $\iota: H\longrightarrow Q^{2k}$ satisfies 
$$\iota(q_1,q_2,\dots,q_k)=(q_1,q_1,q_2,q_2,\dots,q_k,q_k).$$
Since $Q^2/\Delta Q\cong S^3$, we can conclude that $$M\cong Q^{2k}/H\cong (Q^2/\Delta Q)^k\cong (S^3)^k.$$

\paragraph*{\textbf{Part 3: Rigidity for $n=3k+2$.}} The arguments used to analyse this case are analogous to those used for $n=3k$. We have $\Q_{ae}(M)=\lfloor\tfrac{6k+4}{3}\rfloor=2k+1$. Let $Q^{2k+1}$ act almost-effectively on $M$, and hence also almost-effectively on a principal orbit $P$. If $\dim P\leq \dim M-1=3k+1$, then $\Q_{ae}(P)\leq 2k$, which contradicts the existence of an almost-effective $Q^{2k+1}$ action on $P$. Therefore, the codimension of $P$ on $M$ is $0$ and hence $M\cong Q^{2k+1}/H$. Since $Q^{2k+1}$ is simply connected, we have $\pi_1(M)\cong\pi_0(H)$, and since $M$ is simply connected, we can conclude that $H$ is connected and hence $H\cong Q^a\times T^b$ by \cref{connected subgroups of Q^r}. We have
$$3k+2=\dim M= \dim Q^{2k+1}-\dim H=6k+3-3a-b$$ 
and therefore $$3k+1=3a+b.$$

Using the inequality $2k+1\geq 2a+b$, we obtain that $a\geq k$. Since the action is almost-effective, the set $I_i$ satisfy $|I_i|\geq 2$ for all $1\leq i\leq a$ by \cref{connected subgroups of Q^r almost-effective}, which implies that $a\leq k$. Thus, $a=k$ and hence $b=1$. We can conclude that $H\cong Q^k\times S^1$, that $|I_i|=2$ for all $i$ and that $|J|=1$. By \cref{connected subgroups of Q^r}, we can modify the action of $Q^{2k+1}$ by precomposing it with an inner automorphism and a permutation of the factors of $Q^{2k+1}$ so that the sets $I_i$ satisfy that $I_{i}=\{2i-1,2i\}$ for all $1\leq i\leq a$, that $J=\{2k+1\}$, and that the inclusion $\iota: H\longrightarrow Q^{2k+1}$ satisfies 
$$\iota(q_1,q_2,\dots,q_k,z)=(q_1,q_1,q_2,q_2,\dots,q_k,q_k,\eta(z)),$$
where $\eta: S^1\longrightarrow Q$ denotes the inclusion of a maximal torus. Since $Q^2/\Delta Q\cong S^3$ and $Q/S^1\cong S^2$, we can conclude that $$M\cong Q^{2k+1}/H\cong (Q^2/\Delta Q)^k\times (Q/S^1)\cong (S^3)^k\times S^2.$$

\paragraph*{\textbf{Part 4: Rigidity for $n=3k+1$, the homogeneous case.}} Finally, we study the case where $n=3k+1$, hence $\Q_{ae}(M)=\lfloor\tfrac{6k+2}{3}\rfloor=2k$. Let $Q^{2k}$ act almost-effectively on $M$, and hence it also acts almost-effectively on a principal orbit $P$. The main difference between this case and the previous two is that we cannot generally conclude that $Q^{2k}$ acts transitively on $M$. Indeed, if $\dim P\leq \dim M-2=3k-1$, then $\Q_{ae}(P)\leq 2k-1$, which contradicts the existence of an almost-effective $Q^{2k}$ action on $P$. However, if $\dim P=\dim M-1=3k$, then $\Q_{ae}(P)\leq 2k$ and therefore it is possible for the action of $Q^{2k}$ on $M$ to have cohomogeneity one. Thus, we study two separate cases: when the action is transitive and when it has cohomogeneity one.

First, assume the action is transitive, so $M\cong Q^{2k}/H$. Again, since $Q^{2k}$ and $M$ are simply connected, we conclude that $H$ is connected and therefore $H\cong Q^a\times T^b$ by \cref{connected subgroups of Q^r}. We have
$$3k+1=\dim M= \dim Q^{2k}-\dim H=6k-3a-b$$
and therefore
$$3k-1=3a+b.$$
The inequality $2k\geq 2a+b$ implies that $a\geq k-1$. On the other hand, $k\geq a$, so $a=k$ or $a=k-1$. Note that $a=k$ would imply $b=-1$, which is not possible since $b\geq 0$. Thus, $a=k-1$ and $b=2$. Moreover, we have that $H\cong Q^{k-1}\times T^2$, that $|I_i|=2$ for all $i$ and that $|J|=2$. By \cref{connected subgroups of Q^r}, we modify the action of $Q^{2k}$ by precomposing it with an inner automorphism and a permutation of the factors of $Q^{2k}$ so the sets $I_i$ satisfy $I_{i}=\{2i-1,2i\}$ for all $1\leq i\leq k-1$, $J=\{2k-1,2k\}$ and the inclusion of $\iota: H\longrightarrow Q^{2k}$ satisfies
$$\iota(q_1,q_2,\dots,q_{k-1},z_1,z_2)=(q_1,q_1,q_2,q_2,\dots,q_{k-1},q_{k-1},\eta(z_1),\eta(z_2)),$$
where $\eta: S^1\longrightarrow Q$ denotes the inclusion of a maximal torus. We can conclude that $$M\cong Q^{2k}/H\cong (Q^2/\Delta Q)^{k-1}\times (Q/S^1)^2\cong (S^3)^{k-1}\times S^2\times S^2.$$

\paragraph*{\textbf{Part 5: Rigidity for $n=3k+1$, the cohomogeneity one case.}} We assume now that the action of $Q^{2k}$ on $M$ has cohomogeneity one. Since $M$ is simply connected, $M$ has two non-principal orbits. Let $H$ denote the isotropy subgroup of a principal orbit and $K_+$ and $K_-$ denote the isotropy subgroups of the two non-principal orbits. Then $M$ is diffeomorphic to the union of two disk bundles $Q^{2k}\times_{K_+}\mathbb{D}^{l_++1}\cup_{Q^{2k}/H}Q^{2k}\times_{K_-}\mathbb{D}^{l_-+1}$. 

Our first goal is to prove that $l_\pm \neq 0$. The transitive action of $Q^{2k}$ on $Q^{2k}/H$ is almost-effective, and hence it is also almost-effective on the cover $Q^{2k}/H_0$, where $H_0$ denotes the connected component of the identity. Since $\dim Q^{2k}/H_0=3k$, we can conclude that $H_0\cong Q^k$ and after precomposing the action of $Q^{2k}\longrightarrow\Diff(M)$ by an inner automorphism and a permutation of the factors of $Q^{2k}$, we can assume that the inclusion morphism $\iota_0:H_0\longrightarrow Q^{2k}$ satisfies 
$$\iota_0(q_1,q_2,\dots,q_k)=(q_1,q_1,q_2,q_2,\dots,q_k,q_k)$$ 
and, in particular, $Q^{2k}/H_0\cong (S^3)^k$. Since $H_0$ is a normal subgroup of $H$, we can conclude that $H\leq N_{Q^{2k}}(H_0)$. A straightforward computation shows that $N_{Q^{2k}}(H_0)/H_0\cong (\Z/2)^k$ and it is generated by the projection of elements $x_i=(x_{i,1},\dots,x_{i,2k})\in Q^{2k}$ satisfying $x_{i,2i}=-Id$ and $x_{i,j}=Id$ otherwise.

Since $M$ is simply connected, we have $l_\pm\neq 0$ by, \cite[Proposition 1.7]{hoelscher2010classification}. Our next goal is to prove that $l_\pm\neq 1$. To reach a contradiction, assume $l_+=1$, and hence $K_+/H\cong S^1$. Note that $H_0$ is a closed subgroup of $K_{+0}$ and therefore $K_{+0}/H_0\cong S^1$. By \cref{connected subgroups of Q^r}, we have $K_{+0}\cong Q^{a_+}\times T^{b_+}$. Using that $H_0\cong Q^k$ and that $K_{+0}/H_0\cong S^1$, we can conclude that $a_+=k$ and $b_+=1$. Let $I_i^+$ and $J^+$ be the sets defined in \cref{connected subgroups of Q^r}, which characterise the inclusion $\iota_{+0}:K_{+0}\longrightarrow Q^{2k}$ up to inner automorphism. Note that if we precompose $\iota_{+0}$ with the inclusion $H_0\longrightarrow K_{+0}$ we obtain the inclusion $\iota_0:H_0\longrightarrow Q^{2k}$. Since $H_0$ is mapped into the $Q^k$-component of $K_{+0}$, we have that $[2k]=\bigcup_{i=1}^k I_i\subseteq\bigcup_{i=1}^k I_i^+\subseteq[2k]$ and therefore $\bigcup_{i=1}^k I_i^+=[2k]$. Thus, $J^+$ is empty since $J^+=[2k]\setminus\bigcup_{i=1}^kI_i^+$, contradicting the fact that $b_+=1$. We have $l_+ \neq 1$, and the same argument shows that $l_-\neq 1$.

We have seen that $l_\pm\geq 2$, which implies that $\pi_1(M)\cong \pi_1(Q^{2k}/H)\cong\pi_0(H)$. Since $M$ is simply connected, $H$ is connected. Moreover, since $K_\pm/H\cong S^{l_\pm}$ with $l_\pm \geq 2$, the groups $K_\pm$ are also connected and by \cref{connected subgroups of Q^r}, $K_\pm\cong Q^{a_\pm}\times T^{b_\pm}$. Since $K_\pm/H\cong S^{l_\pm}$ with $l_\pm\geq2$, the long exact sequence of homotopy yields $b_\pm=0$. Note that $l_\pm=3a_\pm-3k$ and hence $l_\pm$ are divisible by $3$ and $a_\pm\geq k+1$. If $l_\pm\geq 6$, then the long exact sequence of homotopy would induce an isomorphism 
$$\Z^{a_\pm}\cong\pi_3(K_\pm)\cong \pi_3(H)\cong \Z^{k},$$ 
which is not possible since $a_\pm\neq k$. Hence the only possibility is $l_\pm=3$ and $a_\pm=k+1$. 

Let $\overline{\iota}_\pm: Q^k\longrightarrow Q^{k+1}$ and $\iota_\pm: Q^{k+1}\longrightarrow Q^{2k}$ denote the inclusion morphisms of the isotropy subgroups given by the cohomogeneity one action of $Q^{2k}$ on $M$. Let $\overline{I}_i^{\pm}$ and $I_i^{\pm}$ be the sets defined in \cref{connected subgroups of Q^r} that characterize each of the four inclusions up to inner automorphism. Recall that we can assume without loss of generality that $\iota=\iota_+\circ\overline{\iota}_+=\iota_-\circ\overline{\iota}_-$ satisfies $$\iota(q_1,q_2,\dots,q_k)=(q_1,q_1,q_2,q_2,\dots,q_k,q_k),$$ and thus the sets defined in \cref{connected subgroups of Q^r} are $I_i=\{2i-1,2i\}$. The equalities $\iota=\iota_+\circ\overline{\iota}_+=\iota_-\circ\overline{\iota}_-$ imply that for all $1\leq i\leq k$ we have  
$$I_i=\bigcup_{\overline{i}\in \overline{I}^+_i} I^+_{\overline{i}} =\bigcup_{\overline{i}\in \overline{I}^-_i} I^-_{\overline{i}}.$$

This implies that there exist $\overline i_+, \overline i_-\in[k]$ satisfying that 
$$|\overline I^\pm_j|=\begin{cases}
    2 \text{ if $j=i_\pm$} \\ 1 \text{ otherwise.}
\end{cases} $$

\paragraph*{\textbf{Part 6: Rigidity for $n=3k+1$, case 2.(b).}} From here, we distinguish two cases. First, we assume that $\overline i_+=\overline i_-$. Without loss of generality, assume $i_+=i_-=1$. Indeed, we can change the action of $Q^{2k}$ on $M$ with a permutation $Q^{2k}\longrightarrow Q^{2k}$ which preserves the pairs $\{2i-1,2i\}$ and exchange the pair $\{1,2\}$ with $\{2\overline i_+-1,2\overline i_+\}$. We need to describe explicitly the inclusions $\overline{\iota}_\pm:Q^k\longrightarrow Q^{k+1}$ and $\iota_\pm:Q^{k+1}\longrightarrow Q^{2k}$, taking into account the possible conjugations appearing in \cref{connected subgroups of Q^r}. Thus, there exists $x_\pm=(x_{1,\pm},\overline x_{1,\pm},x_{2,\pm},x_{3,\pm},\dots,x_{k,\pm})\in Q^{k+1}$ and $y_\pm=(y_{1,\pm},\dots,y_{2k,\pm})\in Q^{2k}$ such that $$ \overline \iota_\pm(q_1,\dots,q_k)=C_{x_\pm}(q_1,q_1,q_2,q_3,\dots,q_k) $$ 
and $$\iota_\pm(q_1,\overline q_1,q_2,q_3,\dots,q_k)=C_{y_\pm}(q_1,\overline q_1,q_2,q_2,\dots, q_k,q_k).$$

We note that by \cref{cohomogeneity one equivalent}.2 we can conjugate all the groups by the element $y_+^{-1}$, obtaining a cohomogeneity one manifold $Q^{2k}$-equivariantly diffeomorphic to the original one. Hence we can assume without loss of generality that $y_+$ is trivial, and $x_\pm$ and $y_-$ change accordingly to the above transformation.

The equality $\iota_+\circ \overline \iota_+=\iota_-\circ \overline\iota_-$ implies $$C_{x_{1,+}}(q_1)= C_{y_{1,-}x_{1,-}}(q_1) \qquad \text{and} \qquad C_{\overline x_{1,+}}(q_1)= C_{y_{1,-}\overline x_{1,-}}(q_1)$$ for the first coordinate, and $$ C_{y_{(2j-1),-} x_{j,-}}(q_j)=C_{x_{j,+}}(q_j)=C_{y_{(2j),-} x_{j,-}}(q_j) $$ for all $2\leq j\leq k$. The second equality implies $$y_{(2j-1),-}y_{(2j),-}^{-1}\in\Zc Q.$$

We now apply operation 3 of \cref{cohomogeneity one equivalent} to the negative part, replacing $K_-$ by $hK_-h^{-1}$ where $h=(1,1,y_{3,-}^{-1},y_{3,-}^{-1},\dots,y_{2k-1,-}^{-1},y_{2k-1,-}^{-1})\in H\subseteq N_{Q^{2k}}(H)_0$, while keeping $H$, $K_+$ and $Q^{2k}$ fixed. This gives a $Q^{2k}$-equivariantly diffeomorphic manifold. Using the original coordinates on $K_-$, the new inclusion $\iota_-'$ is $C_h\circ\iota_-$. The new inclusion $\overline \iota_-'$ is the inclusion of the fixed subgroup $H$ into $hK_-h^{-1}$. In coordinates, it is uniquely determined by requiring that its composition with the new $\iota_-'$ inclusion is still $\iota$. Such an inclusion exists because $h$ normalizes $H$, and it is unique because the outer inclusion is injective. We retain the notation $\iota_-$ and $\overline\iota_-$ for the resulting maps. 

Since $y_{(2j-1)-}y_{(2j)-}^{-1}\in\Zc Q$, the new coordinates for $hK_-h^{-1}$ are $$\begin{aligned} C_{y_{(2j-1),-}^{-1}}\bigl(C_{y_{(2j-1),-}}(q_j)\bigr)&=q_j,\\ C_{y_{(2j-1),-}^{-1}}\bigl(C_{y_{(2j),-}}(q_j)\bigr)&=C_{y_{(2j-1),-}^{-1}y_{(2j),-}}(q_j)=q_j, \end{aligned}$$ for all $j\geq2$. Therefore, $$\iota_-(q_1,\overline q_1,q_2,\dots,q_k)=\bigl(C_{y_{1,-}}(q_1),C_{y_{2,-}}(\overline q_1),q_2,q_2,\dots,q_k,q_k\bigr).$$ 
In particular, both inclusions $\iota_+$ and $\iota_-$ are now exactly the standard diagonal inclusion on every common factor $j\geq2$. We also retain the notation $K_-$ for the new group $hK_-h^{-1}$.

We now focus on the inclusions of $\overline \iota_\pm:H\longrightarrow K_\pm$. Since $\iota_+\circ \overline \iota_+=\iota_-\circ \overline\iota_-$, we have $$C_{x_{j,-}}=C_{x_{j,+}},$$ for $j\geq 2$. The same argument as above gives $$x_{j,+}^{-1}x_{j,-}\in\Zc Q.$$ We can precompose $\overline\iota_\pm$ and $\iota$ with $$C_{(1,x_{2,+}^{-1},\dots,x_{k,+}^{-1})}:H\longrightarrow H,$$ since the resulting diagram is isomorphic to the original one via this automorphism on $H$ and the identity maps on $K_\pm$ and $G$, and hence defines a $G$-equivariantly diffeomorphic manifold. For each $j\geq2$, the resulting inclusions are given by $C_{x_{j,+}x_{j,+}^{-1}}=\operatorname{Id}_Q$ and $C_{x_{j,-}x_{j,+}^{-1}}=\operatorname{Id}_Q$, since $x_{j,-}x_{j,+}^{-1}\in\Zc Q$. Thus, the new inclusions $\overline \iota_\pm':H\longrightarrow K_\pm$ satisfy $$(q_1,\dots,q_k)\mapsto (x_{1,\pm}q_1x_{1,\pm}^{-1},\overline x_{1,\pm}q_1\overline x_{1,\pm}^{-1},q_2,q_3,\dots,q_k).$$

Using the identifications already chosen, we decompose $$H=Q^{\{1\}}\times Q^{\{2,\dots,k\}},\qquad K_\pm=Q^2\times Q^{\{2,\dots,k\}},\qquad Q^{2k}=Q^2\times Q^{\{3,\dots,2k\}}.$$ On $Q^{\{2,\dots,k\}}$ both inclusions $H\longrightarrow K_\pm$ are the identity on $Q^{\{2,\dots,k\}}$, and both inclusions $K_\pm\longrightarrow Q^{2k}$ restrict to $$Q^{\{2,\dots,k\}}\longrightarrow Q^{\{3,\dots,2k\}},\qquad(q_2,\dots,q_k)\longmapsto(q_2,q_2,\dots,q_k,q_k).$$ Thus, we can use \cref{cohomogeneity one product} and \cref{cohomogenity one structure} to conclude that $M$ is diffeomorphic to the product of a four-dimensional cohomogeneity one manifold and $$(Q^2/\Delta Q)^{k-1}\cong(S^3)^{k-1}.$$

The action on this four-dimensional factor is determined by the inclusions $$Q\longrightarrow Q^2\longrightarrow Q^2$$ restricted to the $Q^{\{1\}}$ factor. They are given by $$q_1\longmapsto\bigl(C_{x_{1,+}}(q_1),C_{\overline x_{1,+}}(q_1)\bigr),\qquad(\overline q_1,\overline q_2)\longmapsto (\overline q_1,\overline q_2)$$ on the positive part and by $$q_1\longmapsto\bigl(C_{x_{1,-}}(q_1),C_{\overline x_{1,-}}(q_1)\bigr),\qquad(\overline q_1,\overline q_2)\longmapsto\bigl(C_{y_{1,-}}(\overline q_1),C_{y_{2,-}}(\overline q_2)\bigr)$$ on the negative part.

Applying the transformations in \cref{cohomogeneity one equivalent} again, we can modify the action on the four-dimensional factor so that the resulting action has diagram 
$$\Delta Q\subseteq Q^2,Q^2\subseteq Q^2.$$
The inclusions in the positive and negative parts coincide, and they satisfy $\overline \iota_\pm=\Delta$ and $\iota_\pm=id_{Q^2}$. This is the standard action on $S^4\subseteq \mathbb H\oplus\R$, given by $(a,b)\cdot(q,t)=(aqb^{-1},t)$ (see \cite[Proposition 1.20]{hoelscher2010classification}). Consequently, $$M\cong S^4\times(S^3)^{k-1}.$$

\paragraph*{\textbf{Part 7: Rigidity for $n=3k+1$, case 2.(c).}} Assume now that $i_+\neq i_-$. By permuting the factor pairs as in the preceding case, we may assume $i_+=1$ and $i_-=2$. Applying the same arguments as above to the common factors $j\geq3$, we can assume that the inclusions $\overline\iota_\pm:H\longrightarrow K_\pm$ are identity maps on the factor $Q^{\{3,...,k\}}$ and that $\iota_\pm:K_\pm\longrightarrow Q^{2k}$ restricted on the image of $Q^{\{3,...,k\}}$ are the standard diagonal embeddings. 
Thus, we obtain the decompositions $$H=Q^2\times Q^{k-2},\qquad K_\pm=Q^3\times Q^{k-2},\qquad Q^{2k}=Q^4\times Q^{2k-4}.$$ By \cref{cohomogeneity one product}, we have $$M\cong N\times(Q^2/\Delta Q)^{k-2}\cong N\times(S^3)^{k-2},$$ where $N$ is a closed, simply connected seven-dimensional manifold with a cohomogeneity one action of $Q^4$. We obtain its inclusion maps by restricting to the first two factors of $H$ and the first four coordinates of $Q^{2k}$.

We now describe the action on the seven-dimensional factor $N$. After the initial simultaneous conjugation setting $y_+=1$, its positive part is given by $$\overline\iota_+(q_1,q_2)=\bigl(C_{x_{1,+}}(q_1),C_{\overline x_{1,+}}(q_1),C_{x_{2,+}}(q_2)\bigr),\qquad\iota_+(\overline q_1,\overline q_2,\overline q_3)=(\overline q_1,\overline q_2,\overline q_3,\overline q_3),$$ and its negative part by $$\overline\iota_-(q_1,q_2)=\bigl(C_{x_{1,-}}(q_1),C_{x_{2,-}}(q_2),C_{\overline x_{2,-}}(q_2)\bigr),$$ $$\iota_-(\overline q_1,\overline q_2,\overline q_3)=\bigl(C_{y_{1,-}}(\overline q_1),C_{y_{2,-}}(\overline q_1),C_{y_{3,-}}(\overline q_2),C_{y_{4,-}}(\overline q_3)\bigr).$$ 

The equality $\iota_+\circ\overline\iota_+=\iota_-\circ\overline\iota_- $ gives $$C_{y_{1,-}x_{1,-}}=C_{x_{1,+}},\qquad C_{y_{2,-}x_{1,-}}=C_{\overline x_{1,+}},$$ and $$C_{y_{3,-}x_{2,-}}=C_{x_{2,+}}=C_{y_{4,-}\overline x_{2,-}}.$$ Applying operation 2 of \cref{cohomogeneity one equivalent}, we conjugate all groups by $$(x_{1,+}^{-1},\overline x_{1,+}^{-1},x_{2,+}^{-1},x_{2,+}^{-1}).$$  The inclusion of $H$ into $Q^4$ then becomes $$(q_1,q_2)\longmapsto(q_1,q_1,q_2,q_2).$$

We consider now the restricted action of $G=Q^{\{1\}}\times Q^{\{3\}}\subseteq Q^4$. The action of $G$ on a principal orbit $Q^4/(\Delta Q\times \Delta Q)$ is transitive and has trivial isotropy subgroup. Consequently, a principal orbit of the action of $G$ on $N$ has dimension $6$, so the action of $G$ on $N$ is also of cohomogeneity one. The singular orbits of the action of $G$ on $N$ are determined by the inclusions $\iota_\pm$ and $\overline\iota_\pm$, since they coincide with the singular orbits of the action of $Q^4$. At the positive singular orbit, the first pair of coordinates varies independently while the second pair forms a conjugate diagonal subgroup. Thus, the restricted isotropy is $Q\times\{1\}$. Similarly, at the negative singular orbit, the restricted isotropy is $\{1\}\times Q$. Thus, the restricted action has diagram 
$$\{1\}\subseteq Q\times\{1\},\{1\}\times Q\subseteq Q^2.$$
By the classification in \cite[Section 4, Case $1_7$A]{hoelscher2010classification}, this is the standard sum action on $S^7\subseteq\mathbb H\oplus\mathbb H$ (see \cite[1.21. Special types of actions.]{hoelscher2010classification}). Consequently, 
$$M\cong S^7\times(S^3)^{k-2}.$$ 
This completes the proof of \cref{bound almost-effective quaternionic degree of symmetry}.
\end{proof}

\printbibliography

\end{document}